\documentclass[11pt,a4paper,leqno,twoside, headinclude,reqno]{amsart}

\usepackage[tracking=true]{microtype}
\DeclareMicrotypeSet*[tracking]{my}%
  { font = */*/*/sc/* }%
\SetTracking{ encoding = *, shape = sc }{ 45 }

\title[Counter-Examples to a generalised TRC]{Counter-Examples to a generalised\\ Toral Rank Conjecture}
\def\titl{Counter-Examples to a generalised\\ Toral Rank Conjecture}
\def\auth{Manuel Amann}
\date{November 24th, 2020}


\author{\auth}

\usepackage[english]{babel}

\usepackage[colorlinks,pdftex, plainpages=false]{hyperref}
\hypersetup{pdftitle=\titl, pdfauthor=\auth, pdftoolbar=false,
plainpages=false, hyperindex=true, pdfdisplaydoctitle=true}

\hypersetup{colorlinks,%
citecolor=black,%
filecolor=black,%
linkcolor=black,%
urlcolor=black,%
pdftex}

\usepackage{color}

\usepackage{amsmath, amssymb, amscd, amsthm}
\usepackage{stmaryrd}
\usepackage{latexsym}
\usepackage[all]{xy}
\usepackage{pb-diagram}
\usepackage{rotating}
\usepackage{multicol}
\usepackage{lscape}
\usepackage{wasysym}
\usepackage{longtable}
\usepackage{enumerate}
\usepackage{arydshln}

\xyoption{all}

\newtheorem{theo}{Theorem}[section]
\newtheorem{main}{Theorem}

\newtheorem*{main*}{Theorem}

\newtheorem*{mainprop*}{Proposition}
\renewcommand{\themain}{\Alph{main}}

\newtheorem{mainconj}{Conjecture}

\newtheorem{prop}[theo]{Proposition}
\newtheorem{defi2}[theo]{Definition}
\newtheorem*{defi2*}{Definition}

\newenvironment{defi*}{\begin{defi2*}\normalfont}{\end{defi2*}}
\newenvironment{defin}[1]{\begin{defi2}[#1]\normalfont}{\end{defi2}}
\newenvironment{defin*}[1]{\begin{defi2*}[#1]\normalfont}{\end{defi2*}}
\newtheorem*{rem2*}{Remark}
\newenvironment{rem*}{\begin{rem2*}\normalfont}{\hfill$\boxbox$\end{rem2*}}
\newtheorem{rem2}[theo]{Remark}
\newenvironment{rem}{\begin{rem2}\normalfont}{\hfill$\boxbox$\end{rem2}}

\newtheorem{lemma}[theo]{Lemma}

\newtheorem*{cor*}{Corollary}

\newtheorem{conj}[theo]{Conjecture}

\newtheorem*{conj*}{Conjecture}
\newtheorem*{theo*}{Theorem}
\newtheorem*{ques*}{Question}
\newtheorem*{mi2}{Main Idea}

\newtheorem{ex2}[theo]{Example}

\newenvironment{exa}[1]{\begin{ex2}[#1]\normalfont}{\hfill$\boxbox$\end{ex2}}
\newtheorem{exer2}[theo]{Exercise}

\newtheorem{alg2}[theo]{Algorithm}

\newtheorem{conv2}[theo]{Convention}
\newenvironment{conv}{\begin{conv2}\normalfont}{\hfill$\boxbox$\end{conv2}}

\newcommand{\cc}{{\mathbb{C}}}                                     
\newcommand{\nn}{{\mathbb{N}}}                                     
\newcommand{\qq}{{\mathbb{Q}}}                                     
\newcommand{\rr}{{\mathbb{R}}}                                     
\newcommand{\pp}{{\mathbf{P}}}                                     
\newcommand{\s}{{\mathbb{S}}}                                      
\newcommand{\zz}{{\mathbb{Z}}}                                     
\newcommand{\GL}{{\mathbf{GL}}}                                    
\newcommand{\SU}{{\mathbf{SU}}}                                    
\newcommand{\Sp}{{\mathbf{Sp}}}                                    
\newcommand{\B}{{\mathbf{B}}}                                      
\newcommand{\E}{{\mathbf{E}}}                                      
\newcommand{\dif} {{\operatorname{d}}}                             
\newcommand{\In} {{\,\subseteq\,}}                                 
\newcommand{\im} {{\operatorname{im\,}}}                           
\newcommand{\APL}{{\operatorname{A_{PL}}}}                         
\newcommand{\rk}{{\operatorname{rk\,}}}                            
\newcommand{\Der}{{\operatorname{Der}}}                            
\newcommand{\aut}{{\operatorname{aut}}}                            
\newcommand{\co}{\colon\thinspace}                                 
\newcommand{\vproof}{{\begin{flushright} \qed                      
                      \end{flushright}}}

\newcommand{\comment}[1]{}                                         
\newcommand{\xto}[1]{\xrightarrow{#1}}                             
\newcommand{\hto}[1]{\overset{#1}{\hookrightarrow}}                

\newcommand{\ack}{\noindent\textbf{Acknowledgements. }}            
\newcommand{\str}{\noindent\textbf{Structure of the article. }}    
\newcommand{\even}{\textrm{even}}                                  

\newenvironment{prf}{\begin{proof}[\textsc{Proof}]} {\end{proof}}     

\begin{document}

\maketitle \thispagestyle{empty}


\begin{abstract}
The toral rank conjecture speculates that the sum of the Betti numbers of a compact manifold admitting a free action of a torus of rank $r$ is bounded from below by $2^r$. Clearly, such an action yields a torus bundle, and, more generally, the same cohomological bound is conjectured for total spaces of suitable topological torus fibrations by F\'elix--Oprea--Tanr\'e.

In this article we show that this generalised toral rank conjecture cannot hold by providing various different counter-examples to it (for each rank $r\geq 5$). In particular, we show that there are sequences of smooth nilpotent fibre bundles of nilmanifolds with fibre a torus of rank $r$ such that the quotient of the total dimensions of the cohomologies of total space and fibre
even converges to $0$ with $r$ tending to infinity. 

We moreover prove that none of our torus fibrations can be realised by almost free torus actions. More precisely, in the depicted sequence the difference of the ranks of the torus fibres in the bundle and the toral ranks (which are constant one) even tends to infinity.

Similar to recent examples by Walker (of a completely different nature) this shows that the toral rank conjecture is not likely to follow from ``weaker conjectures or structures''.
\end{abstract}


\section*{Introduction}

Symmetry plays an extraordinary role in our understanding of geometry---both does it incorporate simplicity and inherent beauty of an object. In particular, Lie group actions constitute a most evident way of formulating such ``self-similarities''.
All compact connected Lie groups possess maximal tori, and actions of compact tori (most naturally reflecting our notion of turning the object) are amongst the first actions of groups of positive dimension to naturally come into play. Such torus actions play an overly important role. Their theory, behaviour, and geometric applications vary strongly with the issue whether they have only finite isotropy, i.e.~they are \emph{almost free}, or whether they come with higher-dimensional stabiliser groups. In the case of almost free actions the following conjecture originally due to Steve Halperin (see \cite{Hal85}) (and later having been generalised as the Halperin--Carlsson conjecture to finite characteristic) has initiated a large amount of fruitful research programmes. However, till today it still seems open in general.
\begin{conj}[toral rank conjecture]
Let $M$ be a simply-connected compact manifold equipped with an almost free
action of a compact torus $T^r=\s^1\times \stackrel{(r)}{\ldots} \times \s^1$ of rank $r$. Then the sum of all Betti numbers, the total
cohomological dimension, satisfies $\dim H^*(M;\qq)\geq 2^r$.
\end{conj}

Recently, in \cite{IW18} amongst others an example
of a finitely generated free differential graded $R$-module (with $R$ a polynomial ring in
$r$ variables---for us playing the role of $H^*(\B T^r)$) is given the total rank of which is less than $2^r$, but which has finite dimensional
(non-zero) cohomology. In particular, this shows that the toral rank conjecture cannot be proven without multiplicative data (see \cite{RS20} for a respective non-realisability result as topological spaces in finite characteristic).

It is the goal of this article to provide various examples which reveal that the toral rank conjecture (rationally) cannot be proved in the category of arbitrary torus fibrations---even under additional stronger assumptions. If the conjecture is true, then this can only be derived using the fact that torus actions are classified by $\B T$, and not merely by the fact that a free torus action $T\curvearrowright X$ yields a bundle $T\hto{} X\to X/T$. The latter, in contrast, was speculated to hold true in \cite[Conjecture 9.71, Problem 9.72, p.~388]{FOT08} by F\'elix--Oprea--Tanr\'e in their seminal textbook \emph{``Algebraic models in Geometry''}. More precisely, this is stated as
\begin{conj}[generalised toral rank conjecture, F\'elix--Oprea--Tanr\'e,]\label{conj01}
If $F\hto{} E\to B$ is a quasi-nilpotent fibration where $F$ is a torus $T^r$ and the rational cohomology of $B$ is finite dimensional, then $\dim H^*(E;\qq)\geq 2^r$. If $F$ is just rationally elliptic is it true that $\dim H^*(E;\qq)\geq \dim H^*(F;\qq)$?
\end{conj}
Recall that a nilpotent space is \emph{rationally elliptic} if both the total dimension of its rational cohomology and the total dimension of its rational homotopy groups are finite.

This article is devoted to providing a variety of differently constructed counter-examples to this conjecture. In particular, we have counter-examples which satisfy many additional properties.
\begin{itemize}
\item
For example, they are indeed smooth nilpotent fibre bundles of nilmanifolds.
\item
The examples do not require rational coefficients and work with certain $\zz_{(p)}$-coefficients already (where $p$ is sufficiently large, but still comparatively small with respect to the torus rank). For every sufficiently large $p$ there is a counter-example in this series. (By $\zz_{(p)}$ we denote the localisation of $\zz$ at the prime $p$.)
\item There is a sequence of examples such that the ratio between the sum of the Betti numbers of the fibre and the one of the total space even becomes arbitrarily small with the rank of the torus increasing.
\item For every torus rank $r\geq 5$ there is a counter-example.
\end{itemize}

Before providing more ``abstractly constructed'' examples in particular for tori of low rank, we begin with ``geometrically motivated'' examples, which result from considering and ``splitting up'' the nilmanifold of upper triangular matrices modulo its standard integral lattice. We denote by $\mathbb{U}(n,\zz)$ the nilpotent group upper triangular matrices with integral entries and $1$ on the diagonal, and by $\mathbb{U}_k(n,\zz)$ its normal subgroup of matrices (with $1$ on the diagonal) and non-trivial only from the $k$-th off-diagonal only (see Section \ref{sec01})---the analogue is denoted by $\mathbb{U}(n,\rr)$, $\mathbb{U}_k(n,\rr)$ when considering real matrix entries.

\begin{main}\label{theoC}
Let $R\in \{\zz_{(p)},\qq\}$ for some prime $p\geq n-1$. Suppose that $k\geq n/2+1$. Then there is a nilpotent fibre bundle of connected nilpotent rationally elliptic spaces
\begin{align*}
T^{\frac{(n-k+1)(n-k+2)}{2}} \times \rr^{\frac{(2- k) (-1 + k - 2 n)}{2}} \hto{} \mathbf{B} \mathbb{U}(n,\zz) \to \mathbf{B}\big( \mathbb{U}(n,\zz)/\mathbb{U}_k(n,\zz) \big)
\end{align*}
with $\rk_{R} H^*\big(\mathbf{B}\big( \mathbb{U}(n,\zz)/\mathbb{U}_k(n,\zz)\big); R \big)<\infty$ and $T$ a compact smooth torus such that
\begin{align*}
\rk_{R}H^*\big(\mathbf{B} \mathbb{U}(n,\zz);R)\big) < 2^{\frac{(n-k+1)(n-k+2)}{2}}
\end{align*}
whenever $n-k\geq \sqrt{2n \log_2n}$. In particular, both restrictions for $k$ are satisfied for $k=\lceil n/2\rceil +1$ and $n\geq 49$. For this value of $k$ it then even holds that
\begin{align*}
\lim_{n\to \infty}\frac{\dim H^*\big(\mathbf{B} \mathbb{U}(n,\zz);\qq\big)}{\dim H^*\bigg(T^{\frac{(\lceil n/2\rceil)(\lceil n/2\rceil+1)}{2}} \times \rr^{\frac{(1- \lfloor n/2\rfloor) (\lfloor n/2\rfloor - 2 n)}{2}};\qq\bigg)}=0
\end{align*}
\end{main}
An asymptotic version with $\zz_{(p)}$-coefficients would read like: For every $\varepsilon>0$ there exists an $n$ such that for all coefficient rings $\zz_{(p)}$ with $p\geq n-1$, the corresponding quotient of $\zz_{(p)}$-ranks is smaller than $\varepsilon$.

\renewcommand{\themain}{A\textsuperscript{*}}

Considering these examples from a different ``geometric perspective'' one obtains the ``doppelganger result''
\begin{main}\label{theoCC}
Suppose that $k\geq n/2+1$. Then there is a nilpotent fibre bundle of connected nilmanifolds
\begin{align*}
T^{\frac{(n-k+1)(n-k+2)}{2}}  \hto{}E_n \to \mathbf{B}\big( \mathbb{U}(n,\zz)/\mathbb{U}_k(n,\zz) \big)
\end{align*}
with $\dim H^*(E; \qq)<\infty$,
and with $T$ a compact smooth torus such that
\begin{align*}
\dim H^*(E_n;\qq) < 2^{\frac{(n-k+1)(n-k+2)}{2}}
\end{align*}
whenever $n-k\geq \sqrt{2n \log_2n}$. In particular, both restrictions for $k$ are satisfied for $k=\lceil n/2\rceil +1$ and $n\geq 49$. For this value of $k$ it then even holds that
\begin{align*}
\lim_{n\to \infty}\frac{\dim H^*(E_n;\qq\big)}{\dim H^*\big(T^{\frac{(\lceil n/2\rceil)(\lceil n/2\rceil+1)}{2}};\qq\big)}=0
\end{align*}
\end{main}
We actually construct the examples as iterated principal $\s^1$-bundles (see Section \ref{sec05}, cf.~Theorem \ref{theo08}).

\renewcommand{\themain}{\Alph{main}}

\setcounter{main}{1}

\begin{rem*}
Recall that there is a prominent \emph{Lie algebra toral rank conjecture} (see Conjecture \ref{conj02}) which specialises and reformulates the toral rank conjecture in the setting of finite dimensional nilpotent Lie algebras. This conjecture hence is substantially more restrictive and, consequently, better results are known in this context. For example, see \cite{DS88} for a confirmation of the $2$-stage case in this context, and \cite{JL04} for a discussion of the substantially higher difficulties in the not yet confirmed $2$-stage case of the toral rank conjecture.

The examples of Theorems \ref{theoC} and \ref{theoCC} even yield counter-examples for an accordingly \emph{generalised Lie algebra toral rank conjecture} (in the vein of the generalisation of the toral rank conjecture via Conjecture \ref{conj01}). The same holds true for the examples from Theorem \ref{theoA} below.
\end{rem*}
\begin{rem*}
In Theorem \ref{theo02} we generalise these examples to fibrations of arbitrarily highly connected spaces using rational techniques: We first construct a rational model of the fibration, then alter degrees whilst preserving its structure and then realise these constructed algebras and their morphisms as fibrations. Consequently, the newly obtained examples satisfy the exactly same estimates on rational cohomology as the examples from Theorems \ref{theoC} and \ref{theoCC}.
\end{rem*}
\begin{rem*}
We prove that none of these examples is realisable by almost free torus actions (see Convention \ref{conv01} and Remark \ref{rem04} which make this precise and relate it to the actual toral rank conjecture). We show this lack of realisablity in two ways: concretely for the given examples of torus fibrations, and, more generally, we show that the rational toral rank of all these total spaces is just $1$. This implies in particular the stronger statement that any nilpotent compact Hausdorff space (or reasonable space, see Section \ref{subsec03}) of the same rational homotopy type does not admit any almost free torus action for a torus of rank $2$ or larger. Hence the difference between the rank of the torus bundle and the rank of an almost freely acting torus can become arbitrarily large. All this equally holds true for the examples constructed for the first part in the subsequent Theorem \ref{theoA}. We discuss all this at length in Section \ref{subsec03}.
\end{rem*}
We remark further that there are easy counter-examples to Conjecture \ref{conj01} for $\zz_p$-torus fibres instead of $\s^1$-torus fibres and $\zz_p$-cohomology instead of rational coefficients as for example given by a $4$-fold covering space of $\s^1$ (with fibre considered as $\zz_2\oplus \zz_2$).

\bigskip

Applying the geometric realisation method underlying Theorem \ref{theoCC} to a completely different algebraic model, i.e.~constructing iterated principal $\s^1$-bundles over $T^2$, lets us clear the situation for tori of small rank.
\begin{main}\label{theoA}
There is a nilpotent fibre bundle of connected nilmanifolds
\begin{align*}
T^k\hto{} E\to T^2
\end{align*}
with base a compact $2$-torus $T^2=\s^1\times\s^1$ and with fibre a compact $k$-torus $T^k$ with $5\leq k\leq 9$ such that the total space $E$ satisfies $\dim H^*(E;\qq)<\dim H^*(T^k;\qq)$.

Consequently, for every $k\geq 5$ there exists a nilpotent torus bundle with fibre $T^k$ over a torus and with total space $E$ a nilmanifold such that the total space satisfies $\dim H^*(E;\qq)< \dim H^*(T^k;\qq)$.
\end{main}
Clearly, tori are (rationally) elliptic. In Remark \ref{rem01} we explain how to obtain infinitely many arbitrarily highly connected elliptic spaces $F$ of the rational type of a product of $k$ odd spheres in a fibration over a product of two odd-dimensional rational spheres whilst preserving the inequality $\dim H^*(E;\qq)<2^k$ (for $5\leq k\leq 9$) for the respective total spaces $E$.

The result also seems not special to a given fibre torus. We merely start with the smallest rank, $k=5$, choosable in our examples, the estimate probably should hold for all subsequent ranks $k$.

For $k\geq 10$ and the second part of the assertion we draw on respective product bundles with base, fibre and total spaces the respective products. Clearly, the initially computed ranks $5\leq k\leq 9$ are sufficient to provide various further such product examples in all ranks $k\geq 10$, since both fibre and total space cohomology dimensions by the K\"unneth formula are the products of the cohomology dimensions of the respective Cartesian factors.

\vspace{3mm}
Morally, the reason for all the presented  examples to exist, is that on the level of the Sullivan models of arbitrary torus fibrations differentials may be twisted with elements from the base; this is not the case for principal fibrations. Surprising enough, in these situations and contrary to many other ones (for example see the next paragraph), such perturbations seem to make a difference.

We finally remark that not only the toral rank conjecture can be formulated via an estimate on the dimension of cohomology. So given a fibration $F\hto{} E \to B$ of nilpotent spaces, then the cited generalisation of the toral rank conjecture would speculate that $\dim H^*(E)\geq \dim H^*(F)$, which we disproved (actually, clearly only using spaces with vanishing Euler characteristics). The as famous so-called \emph{Halperin conjecture} speculates for a rationally elliptic space $F$ with positive Euler characteristic that the Leray--Serre spectral sequence of any such fibration with fibre $F$ degenerates at the $E_2$-term, or, equivalently, that $\dim H^*(E)=\dim H^*(F)\cdot \dim H^*(B)$. The \emph{Hilali conjecture} speculates that the relation of the total cohomological dimension and the total dimension of the rational homotopy groups of a rationally elliptic space is at least $1$. This illustrates once more how prominently such cohomology estimates feature within Rational Homotopy Theory. In Section \ref{subsec04} we reveal further similarities shared with the Halperin conjecture thereby contrasting them with our counter-examples; then, however, from the point of view of the classifying spaces of the involved fibrations.

\vspace{3mm}

\str In Section \ref{sec00} we recall necessary aspects from rational homotopy theory and transformation groups. Section \ref{sec01} establishes the proof of Theorem \ref{theoC} by constructing the examples via upper triangular matrices. In Section \ref{sec04} we then establish their rational models and use them consequently in order to gain further highly connected examples via rational realisation techniques. Moreover, this sets the ground for the proof of Theorem \ref{theoCC}, which draws upon these rational examples and realises them geometrically in a different way. In Section \ref{sec02} we do likewise for Theorem \ref{theoA} and construct the rational examples underlying it. Section \ref{sec05} then realises these examples using iterated circle bundles and contains the proofs of Theorems \ref{theoCC} and \ref{theoA}. In Section \ref{sec03} we finally discuss that the given examples cannot (even not rationally) be realised by torus actions, and we put the results into further context mainly providing an interpretation via the well-known theory of classifying spaces of fibrations.

\ack The author is grateful to Berrin \c{S}ent\"urk and Marc Stephan for their feedback on a previous version of the article, and to Leopold Zoller for his comments helping to sharpen the presentation of the results.

The author was supported both by a Heisenberg grant and his research grant AM 342/4-1 of the German Research Foundation; he is moreover a member of the DFG Priority Programme 2026.


\section{Preliminaries}\label{sec00}

We provide some background information on the techniques and structures we use. In particular, as far as rational homotopy theory is concerned, this is by no means intended to nor can it provide an introduction to the theory. For the latter we recommend the books and articles \cite{FHT01}, \cite{FOT08}, \cite{FHT15}, and \cite{FH17}.

\subsection{Nilpotence structures}
First, we shall be dealing with \emph{nilpotent groups}, \emph{nilpotent spaces}, \emph{(quasi-)nilpotent fibrations} and \emph{nilmanifolds}.

\begin{defin}{nilpotent group}
A group $G$ is \emph{nilpotent} if its lower central series stabilises, i.e.~if for some $n\in \nn$
\begin{align*}
G_1&=G\unrhd G_2:=[G,G_1]\unrhd G_3:=[G,G_2]\unrhd  \ldots \unrhd  G_n=[G,G_{n-1}]=1
\intertext{or, equivalently, if it has a stabilising upper central series (consequently of the same length)}
1&=C_0\unlhd C_1:=C(G) \unlhd C_2 \unlhd \ldots \unlhd C_n=G
\end{align*}
where $C_{i+1}$ is such that $C_{i+1}/C_i=C(G/C_i)$.
\end{defin}

For the definition of nilpotent spaces and nilpotent fibrations we recur on the usual actions of the fundamental group on covering spaces respectively, more generally, on fibres in fibrations. For a concise depiction see \cite[Remark 1.9, p.~6]{TO97}.
\begin{defin}{nilpotent space}
A path-connected space $X$ is \emph{nilpotent} if its fundamental
group $\pi_1(X)$ is a nilpotent group acting nilpotently (with the usual action) on the higher
homotopy groups $\pi_i(X)$ for $i\geq 2$, i.e.~there is a lower central sequence
\begin{align*}
\pi_i(X)=G_1^i\unrhd G_2^i:=[G,G_1^i]\unrhd \ldots \unrhd G_n^i=[G,G_{n-1}^i]=1
\end{align*}
such that the induced $\pi_1(X)$-action on the quotients $G_k^i/G_{k+1}^i$ is trivial.
\end{defin}

The following notions can be found in \cite{Dwy74}.
\begin{defin}{(quasi-)nilpotent fibration}
A fibration of path-connected spaces $F\hto{} E\to B$ is called \emph{nilpotent} if $\pi_1(B)$ acts nilpotently (in the usual way) on $\pi_*(F)$.

A fibration of path-connected spaces is \emph{quasi-nilpotent} if the action of $\pi_1(B)$ on $H^*(F;\zz)$ is nilpotent.
\end{defin}
Recall that a nilpotent space $X$ is \emph{rationally elliptic} if $\dim H^*(X;\qq)<\infty$ and $\sum_{i\geq 2}\dim \pi_{i}\otimes \qq<\infty$.

\begin{rem}\label{rem03}
Let us make some well-known remarks on nilpotent groups respectively some comments circling around the relation between (quasi-)nilpotent fibrations and the involved spaces.
\begin{enumerate}
\item Every subgroup and every quotient of a nilpotent group is nilpotent.
\item A fibration of connected spaces is nilpotent if and only if it is quasi-nilpotent and the fibre is nilpotent (see \cite[Corollary 2.21, p.~71]{Hil76}, cf.~\cite[Lemma 2.18, p.~70]{Hil75}).
\item A connected fibre of a fibration with total space $E$ nilpotent is nilpotent as well (see \cite[Proposition $8^\textrm{bis}$.21, p.~347]{Mcc01}, cf.~\cite[Theorem 2.2, p.~62]{Hil75}).
\item If base and total space are connected nilpotent (CW-complexes) and the fibre is connected, then the fibration is nilpotent (see \cite[Proposition 2.13, p.~67]{Hil75}, cf.~\cite[Remark 2.62, p.~79]{FOT08}).
\end{enumerate}
\end{rem}

Let $F\hto{} E\xto{p} B$  be a (Serre) fibration of path-connected spaces. Let $(\Lambda W,\dif)\xto{m_B} \APL(B)$ be a minimal Sullivan model of $B$, and let $(\Lambda W \otimes \Lambda V,\dif)\xto{m_E} \APL(E)$ be a relative minimal Sullivan model of $\APL(p)\circ m_B$. We call the commutative diagram
\begin{align*}
\xymatrix{
(\Lambda W,\dif) \ar[r]  \ar[d]^{m_B}_\simeq
& (\Lambda W \otimes \Lambda V,\dif) \ar[r]  \ar[d]^{m_E}_\simeq
& (\Lambda V,\bar \dif)  \ar[d]^{m_F}\\
\APL(B)
\ar[r]^{\APL(p)}& \APL(E)
\ar[r]
& \APL(F)
}
\end{align*}
the \emph{minimal model of the fibration}. (Here $\bar \dif$ is induced from $\dif$ by projecting away $\Lambda W$.)

From the more general statements \cite[Theorem 12.1, plus the Remark, p.~25]{FH17} or \cite[Theorem 4.4, p.~30]{TO97} (for quasi-nilpotent fibrations) or from \cite[Theorem 5.1, p.~145]{FHT15} (for ``locally nilpotent'' actions on cohomology) respectively from \cite[Theorem 4.1, p.~445]{GH87} (for nilpotent fibrations) we cite that ``the fibre of the model is the model of the fibre'' under the following circumstances.
\begin{theo}\label{theo03}
In the depicted situation and supposing that the (Serre) fibration is (quasi-)nilpotent with $F$ the homotopy type of a CW-complex (or one of $H^*(F;\qq)$ and $H^*(B;\qq)$ of finite type), we obtain that $m_F$ is a minimal Sullivan model, i.e.~a quasi-isomorphism.
\end{theo}

\bigskip

Let us now introduce some more geometric concepts.
\begin{defin}{nilmanifold}
A \emph{nilmanifold} $N=G/\Gamma$ is the quotient of a simply-connected nilpotent Lie $G$ group by a co-compact lattice $\Gamma$.
\end{defin}
Clearly, if $\Gamma$ is a nilpotent group and $G$ is contractible, then $G/\Gamma\simeq \B \Gamma$ (with fundamental group $\Gamma$ and trivial higher homotopy groups) is a nilpotent space. This will be valid in the situations we are interested in, i.e.~especially for the upper triangular matrices $\Gamma=\mathbb{U}(n,\zz)\In \mathbb{U}(n,\rr)=G$ which we shall quickly recall in the following

\begin{exa}{upper triangular matrices}\label{exa01}
Let $\GL(n,\zz)$ denote the general linear group of all invertible $(n\times n)$-matrices with coefficients in $\zz$. In there we have the subgroup
\begin{align*}
\mathbb{U}(n,\zz)=\left\{
\begin{pmatrix}
1      & x_{2,1} & x_{3,1}   &\ldots & x_{n,1} \\
0      & 1      & x_{3,2}   &\ldots & x_{n,2} \\
\vdots & \ddots      & \ddots  &  \ddots     &  \vdots \\
0      & \ldots &   0     &1      & x_{n,n-1} \\
0      & \ldots &   0     &0      & 1
\end{pmatrix}
\bigg | x_{i,i}\in \zz, 1\leq i\leq n \right\} \In \GL(n,\zz)
\end{align*}
Similarly, we have the nilpotent Lie group $\mathbb{U}(n,\rr)\In \GL(n,\rr)$ of such upper triangular matrices with coefficients in $\rr$ containing $\mathbb{U}(n,\zz)$ as a lattice  and yielding the compact quotient nilmanifold
\begin{align*}
N:=\mathbb{U}(n,\rr)/\mathbb{U}(n,\zz)
\end{align*}
(see \cite[Example 3.17, p.~118]{FOT08}, cf.~\cite{Dwy85}). (Clearly, $\mathbb{U}(n,\zz)$ indeed is a subgroup due to the classical formula for inverses $A^{-1}=\frac{\operatorname{adj} A}{\det A}$, where $\operatorname{adj} A$ is the transpose of the matrix with entries $\tilde a_{i,j}=(-1)^{i+j}\cdot \det A_{i,j}$ where $A_{i,j}$ results from $A$ by cancelling the $i$-th line and the $j$-th column.)

Next consider the subgroup $\mathbb{U}_k(n,\zz)\In \mathbb{U}(n,\zz)$,
\begin{align*}
\left\{
\begin{pmatrix}
1      & 0&\ldots&0& x_{k,1} & x_{k+1,1}   &\ldots & x_{n,1} \\
0      & 1      & 0&\ldots&0& x_{k+1,2}   &\ldots & x_{n,2} \\
\vdots & \ddots      & \ddots&\ddots&&  \ddots  &  \ddots     &  \vdots \\
0      & \ldots &   0     &1      & 0&\ldots&0& x_{n,n-k+1} \\
\vdots      &  &     &\ddots      & 1 &\ddots&  &0 \\
\vdots      &  &      &     & \ddots &\ddots& \ddots &\vdots \\
\vdots      &  &      &     &  &\ddots& \ddots &0 \\
0      & \ldots &   \ldots     &  \ldots   & \ldots & \ldots & 0& 1
\end{pmatrix}
\bigg | x_{i,i}\in \zz, 1\leq i\leq n \right\}
\end{align*}
of all those matrices with $k-2$ upper off-diagonals equal to zero, respectively the analog with coefficients in $\rr$. (Clearly, $\mathbb{U}_2(n,\zz)=\mathbb{U}(n,\zz)$, etc.)
\end{exa}

A model for a nilmanifold can be constructed in different ways. Essentially using that a nilmanifold is a coformal space (i.e.~it has a Lie model with vanishing differential, see \cite[Example 7, p.~334]{FHT01}) and that the Lie bracket is dual to the differential of a minimal Sullivan model, one ends up with (see \cite[Theorem 3.18, p.~118]{FOT08})
\begin{theo}\label{theo05}
If $N/\Gamma$ is a nilmanifold, then the complex $(\Lambda \mathfrak{n}^*,\dif)$ associated to dual of the Lie algebra $\mathfrak{n}$ of $N$ is isomorphic to the minimal model of $N/\Gamma$.
\end{theo}
As for defining the differential $\dif$, we choose a homogeneous basis $(X_i)$ of $\mathfrak{n}$ with dual basis $(x_i)$ of $1$-forms.

The basis of the $X_i$ is ordered by the order of central series extensions. That is, we consider the upper central series
\begin{align*}
1&=C_0\unlhd C_1:=C(N) \unlhd C_2 \unlhd \ldots \unlhd C_n=N
\end{align*}
take the quotient by $C_1$
\begin{align*}
1&= C_1/C_1 \unlhd C_2/C_n \unlhd \ldots \unlhd C_n/C_1
\end{align*}
and consider the central extension
\begin{align*}
C_1\hto{} N\to N/C_1
\end{align*}
(cf.~\cite[p.~122]{FOT08}). (We actually draw on the \emph{refined} upper central series such that extensions are by $\zz$.) The quotient of the central series is the central series of the quotient $N/C_1$, which allows for an inductive approach. This gives a sequence of extensions, and we order the $X_i$ according to their place therein, adding the centre last.

We then make $\dif$ dual to the Lie bracket by $\dif x_k(X_i,X_j)=-x_k([X_i,X_j])$. We extend $\dif$ as a graded derivation. Since $[X_i,X_j]=\sum c_{i,j}^l X_l$ (with the actually rational structure constants $c_{i,j}^l$), we obtain that $\dif x_k(X_i,X_j)=-c_{i,j}^k$ and
\begin{align}\label{eqn12}
\dif x_k=- \sum_{i<j} c_{i,j}^k x_i \wedge x_j
\end{align}
(for the ordered basis $(x_i)_i$).

We remark that although the procedure underlying the construction yields a minimal model by relying on comparison arguments with smooth differential forms, it turns out that the construction actually does yield a minimal Sullivan model over the rationals (see \cite[Remark, p.~68]{Has89}).

\subsection{Almost free torus action}\label{subsec03}

We recall the \emph{Borel fibration} $X\hto{}X_{T} \to \B T$ for an action of the compact torus $T^r$ on the space $X$. The total space $X_T$ is the \emph{Borel construction}, $X_T=H\times_T \B T$, its cohomology $H^*_T(X):=H^*(X_T)$ is the \emph{equivariant cohomology} of the action. We shall essentially use it in the following. Recall that due to its (intrinsic) formality (see \cite[Section 2.7, p.~92]{FOT08}) a model of $\B T$ is $H^*(\B T)=\qq[t_1,\ldots, t_r]$ with $\deg t_i=2$ for $1\leq i\leq r$. With a minimal Sullivan model $(\Lambda V,\dif)$ of $X$ a model of the Borel fibration hence is given as the relative Sullivan algebra
\begin{align}\label{eqn11}
(\Lambda V \otimes \qq[t_1,\ldots, t_r] ,\dif)
\end{align}
over the differential subalgebra $(\qq[t_1,\ldots, t_r],0)$.

Recall that a group action on a topological space $X$ is called \emph{almost free} if all isotropy groups are finite. The \emph{toral rank} $\rk X$ of $X$ is the largest rank of a torus acting almost freely on $X$.

We widen the definition of \emph{rational toral rank} slightly compared to the one in \cite[Definition 7.12, p.~278]{FOT08}. Instead of with finite CW-complexes we will work with what Halperin called \emph{reasonable spaces} (see \cite[p.~293]{Hal85}), whilst adapting the definition to our purposes. That is, such a space $X$ is either \emph{simply-connected} and
\begin{itemize}
\item connected paracompact Hausdorff,
\item locally path-connected and semi-locally simply connected,
\item is of finite cohomological dimension over $\qq$ in the sheaf theoretic sense,  $\dim H^*(X;\qq)<\infty$, and
\item $\varinjlim_{x\in U} H^*(U;\qq)=\qq$ (for the limit over all neighbourhoods $U$ of $x$)
\end{itemize}
or, if $\pi_1(X)\neq 0$, we take $X$ to be compact Hausdorff.

The class of reasonable spaces then comprises all finite CW-complexes.

Connected spaces admit minimal Sullivan models. We cite Hsiang's theorem for reasonable spaces from \cite[Proposition 4.1.7, p.~257]{AP93} and \cite[Proposition 7.17, p.~281]{FOT08} in the non-simply-connected case (cf.~\cite[Proposition 1, p.~69]{AH78}) respectively from \cite[Proposition 4.2, p.~298]{Hal85}.
\begin{theo}[Hsiang]\label{theo04}
Let $X$ be a reasonable space $X$ admitting an action by an $r$-torus $T$. The the action is almost free if and only if $H_T^*(X)<\infty$.

Conversely, given a relative algebra as in \eqref{eqn11}, there exists a reasonable space $X$ with model $(\Lambda V,\dif)$ and an almost free $T^r$-action on it if and only if
\begin{align*}
\dim H(\Lambda V \otimes \qq[t_1,\ldots, t_r],\dif)<\infty
\end{align*}
(In this case $X$ may be taken to be a simply-connected respectively nilpotent finite CW-complex with a free $T$-action.)
\end{theo}
The reasonable space is the ``corresponding one'', i.e.~it is taken in the respective category: it is compact Hausdorff if $V^1\neq  0$ and of the more complicated definition otherwise.

Hence we can make
\begin{defin}{rational toral rank}
The \emph{rational toral rank} $\rk_0 X$ of a space $X$ is the largest rank of a torus acting almost freely on a corresponding reasonable space with minimal model isomorphic to the one of $X$.
\end{defin}
By construction, this value $\rk_0 X$ is actually a rational homotopy invariant, which, due to Hsiang, can be computed via minimal models.

Theorem \ref{theo04} motivates an algebraic transcription of the toral rank conjecture. Algebraically encoding the Borel fibration of an almost free torus action via its minimal model, one speculates the following.
\begin{conj}[Algebraic TRC, see Conjecture 7.2, p.~283, in \cite{FOT08}]\label{conj03}
Let $(\Lambda V,\dif)$ be a minimal Sullivan algebra, and let
\begin{align*}
(\Lambda \langle t_1,\ldots, t_r\rangle,0) \hto{} (\Lambda \langle t_1, \ldots, t_r\rangle \otimes \Lambda V,\dif) \to (\Lambda V,\bar \dif)
\end{align*}
be a relative minimal model with $\deg t_i=2$ for $1\leq i\leq r$, and such that the cohomology algebra $H(\Lambda \langle t_1,\ldots, t_r\rangle \otimes \Lambda V, \dif)$ is finite dimensional. Then $\dim H(\Lambda V,\dif)\geq 2^r$.
\end{conj}
Hence an almost-free torus action $T^r\curvearrowright X$ on a reasonable space $X$ satisfies the prerequisites of the algebraic toral rank conjecture. Conversely, any such relative Sullivan algebra produces a free action on a finite CW-complex $X$ with minimal model $(\Lambda V,\bar \dif)$. For this see \cite[Proposition 7.17, p.~280]{FOT08}, which even yields that (if dimension is not divisible by four, which poses more obstructions) the relative minimal model from the algebraic toral rank conjecture can be realised by a smooth torus action on a compact closed manifold. Hence, either starting from the geometric or the algebraic point of view, at the end both the toral rank conjecture (in this generalised form for reasonable spaces) and the algebraic toral rank conjecture require $\dim H^*(X)\geq 2^r$.

\bigskip

Let us now see that in the special case of nilmanifolds the rational toral rank can be determined easily (see \cite[Theorem 7.28, p.~288]{FOT08}).
\begin{theo}\label{theo06}
Let $M$ be a nilmanifold with associated Lie algebra $L$. Then the rational toral rank equals the center of the Lie algebra, i.e.~$\rk_0(M)=\dim C(L)$.
\end{theo}
Since the rational toral rank, $\rk_0 M$, is a rational invariant, this equality remains true for any space of the same rational type as $M$.
In view of this theorem one may reformulate the toral rank conjecture for nilmanifolds/nilpotent Lie algebras (using Lie algebra cohomology)---see \cite[Conjecture 7.29, p.~289]{FOT08}.
\begin{conj}\label{conj02}
If $L$ is a finite dimensional nilpotent Lie algebra defined over the rational numbers, then
\begin{align*}
\dim H^*(L;\qq)\geq 2^{\dim C(L)}
\end{align*}
\end{conj}

Theorem \ref{theo06} has a weaker but more general counterpart. Let $L:=(\pi_*(\Omega X)\otimes \qq,[\cdot, \cdot ])$ be the homotopy Lie algebra of a topological space $X$. From \cite[Proposition 7.27, p.~287]{FOT08} respectively \cite[Corollary 3.5, p.~9]{AP86} we derive the next theorem using that the proof of \cite[Proposition 7.27, p.~287]{FOT08} relies on algebraic computations with models only once Hsiang's theorem holds in the class of considered spaces.
\begin{theo}\label{theo07}
Let $X$ be a nilpotent reasonable space with $\pi_\even(X)\otimes \qq=0$. Then we obtain that
\begin{align*}
\rk_0(X)\leq \dim C(L)
\end{align*}
\end{theo}
Recall that in the case of a nilmanifold with associated Lie algebra $L$, this Lie algebra $L$ is also the rational homotopy Lie algebra (see \cite[p.~288]{FOT08})---this connects the theorems.


\section{Proof of Theorem \ref{theoC}}\label{sec01}

The examples in Theorem \ref{theoC} will be constructed using the nilmanifold of upper triangular matrices
\begin{align*}
\mathbb{U}(n,\rr)/\mathbb{U}(n,\zz)
\end{align*}
we introduced in Example \ref{exa01}. Hence we need to elaborate on this example further.

Note that
\begin{align*}
\dim \mathbb{U}_k(n,\rr)&=(n-k+1)+(n-k)+(n-k-1)+\ldots +1
\\&= \frac{(n-k+1)(n-k+2)}{2}
\\&=:d(n,k)
\end{align*}

The group $\mathbb{U}_k(n,\rr)$ acts freely on $\mathbb{U}(n,\rr)$ by right-multiplication via only altering the subgroup $\mathbb{U}_k(n,\rr)$; and the group quotient is identical to the vector space quotient yielding the obvious diffeomorphism to $\rr^{\frac{n(n-1)}{2}-d(n,c)}$, i.e.~
\begin{align}\label{eqn06}
\mathbb{U}(n,\rr)/\mathbb{U}_k(n,\rr)\cong\rr^{c(n,k)}
\end{align}
for
\begin{align*}
c(n,k):=&\frac{n(n-1)}{2}-d(n,c)
\\=&\frac{n(n-1)}{2}-\frac{(n-k+1)(n-k+2)}{2}
\\=&\frac{(2- k) (-1 + k - 2 n)}{2}
\end{align*}

Indeed, let us quickly check that $\mathbb{U}_k(n,\rr)$ is a group. It is easy to see that the set $\mathbb{U}_k(n,\rr)$ is multiplicatively closed. By passing to the Lie algebra of $\mathbb{U}(n,\rr)$, which, due to the existence of the matrix logarithm in this special case, via the exponential map is bijective to the group itself, one observes that the preimage of $\mathbb{U}_k(n,\rr)$ are the corresponding matrices with $0$ on the diagonal instead of $1$. Hence inverses are given by the negatives and thus are clearly of the same form again. The arguments restrict to the case of coefficients in $\zz$. (The exponential map is a converging power series on upper triangular matrices.)

In the case $k\geq n/2+1$, since it is easily checked that
\begin{align}\label{eqn03}
[A,B]=AB-BA \in \mathfrak{u}_{k_1+k_2-1}(n,\rr) \\ \nonumber \quad \textrm{for } A\in \mathfrak{u}_{k_1}(n,\rr) \textrm{ and } B\in \mathfrak{u}_{k_2}(n,\rr)
\end{align}
it follows that the Lie algebra $\mathfrak{u}_k(n,\rr)$ of $\mathbb{U}_{k}(n,\rr)$ is abelian. Hence the matrix exponential map satisfies
\begin{align*}
\exp([A,B])=\exp(A)\cdot \exp(B)
\end{align*}
As a consequence we have that matrix multiplication on $\mathbb{U}_k(n,\zz)$ becomes addition of corresponding matrix entries, i.e.~there is the isomorphism of groups
\begin{align*}
(\zz^{\oplus d(n,k)},+)\cong (\mathbb{U}_k(n,\zz),\cdot)
\end{align*}
The analogue
\begin{align*}
(\rr^{\oplus d(n,k)},+)\cong (\mathbb{U}_k(n,\rr),\cdot)
\end{align*}
holds true and, consequently, this induces an isomorphism of groups
\begin{align}\label{eqn04}
(\mathbb{U}_k(n,\rr)/\mathbb{U}_k(n,\zz),\cdot)\cong (\rr^{\oplus d(n,k)}/\zz^{\oplus d(n,k)},+)\cong (T^{\oplus d(n,k)},\cdot)
\end{align}

Hence it remains to prove normality in
\begin{lemma}\label{lemma05}
Suppose that $k\geq n/2+1$.
The matrices
\begin{align*}
\zz^{\oplus d(n,k)}\cong \mathbb{U}_k(n,\zz)\unlhd \mathbb{U}(n,\zz)
\end{align*}
form a normal subgroup isomorphic to $\pi_1\big(T^{d(n,k)}\big)$.
\end{lemma}
\begin{prf}
One either observes the normal subgroup property directly for the matrices, or one proceeds as follows: We pass to the Lie algebra $\mathfrak{u}_k(n,\rr)$ and to the image of $\mathbb{U}_k(n,\zz)$ under the matrix logarithm. The Baker--Campbell--Hausdorff formula establishes that the matrix exponential is a homomorphism up to Lie brackets, i.e.~
\begin{align*}
\exp(A)\cdot \exp(B)=\exp(A+B+R)
\end{align*}
with $R$ a term in Lie brackets of two such integral $A, B \in \mathfrak{u}_k(n,\rr)$ of word-length at least $2$. (The formula yields actual converging and not only formal power series in our case.)

It follows for $\exp(B)\in \mathbb{U}_{k'}(n,\zz)$ with $k'\geq 2$ and $\exp(A)\in \mathbb{U}_{k}(n,\zz)$  with $k,k'\geq 2$ that
\begin{align*}
\exp(B)\cdot \exp(A) \cdot \exp(B^{-1})&=\exp(B + A -B + R)
\\&=\exp(A+R)
\\&\in \mathbb{U}_k(n,\zz)
\end{align*}
with $R\in \mathfrak{u}_{k+k'-1}(n,\rr)\In \mathfrak{u}_{k}(n,\rr)$ due to Property \eqref{eqn03}. The surjectivity of $\exp$ in our case yields the normal subgroup property.
\end{prf}
In the folowing we shall use both notations $\mathbb{U}_k(n,\zz)\cong \zz^{\oplus d(n,k)}\cong \pi_1\big(T^{d(n,k)}\big)$ likewise. However, in a slightly more abstract context we use the fundamental group, as subgroups of the triangular matrices we shall prefer the other notations. With the very same arguments we see that also
\begin{align*}
\mathbb{U}_k(n,\rr)\unlhd \mathbb{U}(n,\rr)
\end{align*}
is a normal subgroup. This lets us prove

\begin{lemma}\label{lemma06}
For $k\geq n/2+1$ there is a diffeomorphism
\begin{align*}
\frac{\mathbb{U}(n,\rr)}{\mathbb{U}_k(n,\zz)}\cong T^{d(n,k)} \times \rr^{c(n,k)}
\end{align*}
identifying the homogeneous space with the normal bundle of \linebreak[4]$T^{d(n,k)}\cong\frac{\mathbb{U}_k(n,\rr)}{\mathbb{U}_k(n,\zz)}\In \frac{\mathbb{U}(n,\rr)}{\mathbb{U}_k(n,\zz)}$.
\end{lemma}
\begin{prf}
Since $\mathbb{U}_k(n,\rr)\unlhd \mathbb{U}(n,\rr)$, we obtain that
\begin{align*}
\frac{\mathbb{U}(n,\rr)}{\mathbb{U}_k(n,\zz)} \cong \frac{\mathbb{U}(n,\rr)}{\mathbb{U}_k(n,\rr)} \cdot \frac{\mathbb{U}_k(n,\rr)}{\mathbb{U}_k(n,\zz)}
\cong \mathbb{U}_k(n,\rr)\times_{\mathbb{U}_k(n,\zz)}(\mathbb{U}(n,\rr)/{\mathbb{U}_k(n,\rr)})
\end{align*}
is the normal bundle over
\begin{align*}
\frac{\mathbb{U}_k(n,\rr)}{\mathbb{U}_k(n,\zz)}\cong  T^{d(n,k)}
\end{align*}
(see Property \eqref{eqn04}) with fibre
\begin{align*}
\mathbb{U}(n,\rr)/\mathbb{U}_k(n,\rr)\cong \rr^{d(n,k)}
\end{align*}
(see the observation afore Equation \eqref{eqn06}). Since the action of $\mathbb{U}_k(n,\zz)$ on $\frac{\mathbb{U}(n,\rr)}{\mathbb{U}_k(n,\rr)}$ is trivial, this yields the triviality of the $\rr^{c(n,k)}$-bundle.
\end{prf}

\begin{lemma}\label{lemma01}
Suppose that $k\geq n/2+1$. Then there is a fibre bundle of connected nilpotent rationally elliptic spaces
\begin{align}\label{eqn02}
T^{d(n,k)} \times \rr^{c(n,k)} \hto{} \mathbf{B} \mathbb{U}(n,\zz) \to \mathbf{B}\big( \mathbb{U}(n,\zz)/\mathbb{U}_k(n,\zz) \big)
\end{align}
(and the fibration then necessarily is nilpotent as well).
\end{lemma}
\begin{prf}
From \cite[p.~348]{Die08} we recall that the exact sequence of groups (see Lemma \ref{lemma05})
\begin{align*}
0\to \pi_1(T^{d(n,k)}) \to  \mathbb{U}(n,\zz) \to  \mathbb{U}(n,\zz)/\pi_1(T^{d(n,k)}) \to 0
\end{align*}
induces a fibre bundle with structure group $ \mathbb{U}(n,\zz)/\pi_1(T^{d(n,k)})$ of associated classifying spaces
\begin{align*}
\mathbf{B}\pi_1(T^{d(n,k)})\hto{} \mathbf{B} \mathbb{U}(n,\zz) \to \mathbf{B}\big( \mathbb{U}(n,\zz)/\mathbb{U}_k(n,\zz) \big)
\end{align*}
constructed as
\begin{align*}
\frac{\mathbf{E}\big( \mathbb{U}(n,\zz)\big)}{\zz^{\oplus d(n,k)}}
\hto{}
\frac{\mathbf{E}\big( \mathbb{U}(n,\zz)/\zz^{\oplus d(n,k)} \big) \times \mathbf{E}  \mathbb{U}(n,\zz)}{\mathbb{U}(n,\zz)}
\to
\frac{\mathbf{E}\big( \mathbb{U}(n,\zz)/\zz^{\oplus d(n,k)} \big)}{ \mathbb{U}(n,\zz)/\zz^{\oplus d(n,k)}}
\end{align*}

Now use
\begin{align*}
\mathbb{U}(n,\rr) =\mathbf{E}\big( \mathbb{U}(n,\zz)\big)
\end{align*}
as a model for the contractible classifying spaces together with the right multiplication action of $\zz^{\oplus \frac{k(k+1)}{2}}=\mathbb{U}_k(n,\zz)$.

This yields the fibre bundle
\begin{align*}
\frac{\mathbb{U}(n,\rr)}{\mathbb{U}_k(n,\zz)}\hto{} \mathbf{B} \mathbb{U}(n,\zz) \to \mathbf{B}\big( \mathbb{U}(n,\zz)/\mathbb{U}_k(n,\zz) \big)
\end{align*}

Using Lemma \ref{lemma06} this becomes
\begin{align*}
 T^{d(n,k)} \times \rr^{c(n,k)}\hto{} \mathbf{B} \mathbb{U}(n,\zz) \to \mathbf{B}\big( \mathbb{U}(n,\zz)/\zz^{\oplus d(n,k)} \big)
\end{align*}

The total space is nilpotent, as it is the nilmanifold $\mathbb{U}(n,\rr)/\mathbb{U}(n,\zz)$. As for the base space we recall that the quotient of a nilpotent group is again nilpotent (see \cite[Corollary $8^{\textrm{bis}}$.19, p.~346]{Mcc01}). The base space hence is the classifying space of a nilpotent group and a $K(\mathbb{U}(n,\zz)/\zz^{\oplus d(n,k)},1)$. Hence it is nilpotent as well, and so is the bundle/fibration due to \cite[Remark 2.62, p.~79]{FOT08}. (Trivially, the fibre is a homotopy torus and nilpotent, consequently.)

All spaces have finite dimensional rational homotopy concentrated in degree $1$ and are rationally elliptic, consequently.
\end{prf}

By $\zz_{(p)}$ we denote localisation of $\zz$ at the prime $p$. The group $\mathbb{U}(n,\zz)$ is built from a torus by successive central extensions (corresponding to adding elements from higher off-diagonals---which form normal $\zz$-subgroups---one after the other). Such an extension sequence on the level of classifying spaces corresponds to forming successive spherical $(\B\zz=\s^1)$-fibrations. Each such extension is of the form $0\to \zz\to E_l\to E_{l-1}\to 0$ where $E_{l-1}$ can be considered as the quotient $E_{l}/\zz$ as well as the corresponding quotient of the entire group $\mathbb{U}(n,\zz)$. In particular, this implies that $\mathbf{B}\big( \mathbb{U}(n,\zz)/\mathbb{U}_k(n,\zz) \big)$ (which is some $E_{l-1}$) has the homotopy type of the total space of finitely many successive $\s^1$-fibrations over a torus, i.e.~of a finite CW-complex. Hence its (co)homology is finitely generated. Rational cohomology is the rationalisation of $\zz_{(p)}$-cohomology (localisation is (left-)exact), and the free $\zz_{(p)}$-summand rationalises to the $\qq$-vector space; this finally yields Assertion (4) in the subsequent Lemma \ref{lemma04}.

The next lemma is instrumental to the proof of Theorem \ref{theoC}.
\begin{lemma}\label{lemma04}
For a prime $p\geq n-1$ it holds that
\begin{enumerate}
\item
\begin{align*}
\rk_{\zz_{(p)}} H^*\big(\mathbf{B} \mathbb{U}(n,\zz);\zz_{(p)}\big)&=\dim_\qq H^*\big(\mathbf{B} \mathbb{U}(n,\zz);\qq)\big) =n!
\\\rk_{\zz_{(p)}} H^*(T^{d(n,k)};\zz_{(p)})&=\dim_\qq H^*( T^{d(n,k)};\qq)=2^{d(n,k)}
\end{align*}
and, consequently,
\item
\begin{align}\label{eqn05}
\rk_{\zz_{(p)}}H^*\big(\mathbf{B} \mathbb{U}(n,\zz);\zz_{(p)}\big) < \rk_{\zz_{(p)}} H^*(T^{d(n,k)} \times \rr^{c(n,k)};\zz_{(p)})
\end{align}
for $n-k\geq \sqrt{2n \log_2n}$ and $k\geq n/2+1$; respectively the same with $\qq$-coefficients. In particular, this is satisfied for $k=\lfloor n/2\rfloor +1$ and $n\geq 49$.
\item
It then even holds that
\begin{align}\label{eqn07}
\lim_{n\to \infty}\frac{\dim H^*\big(\mathbf{B} \mathbb{U}(n,\zz);\qq\big)}{\dim H^*(T^{d(n,\lfloor n/2\rfloor +1)} \times \rr^{c(n,\lfloor n/2\rfloor +1)};\qq)}=0
\end{align}
\item
Moreover, we have that
\begin{align*}
\infty&>\rk_{\zz_{(p)}} H^*\big(\mathbf{B}\big( \mathbb{U}(n,\zz)/\zz^{\oplus d(n,k)} \big);\zz_{(p)}\big)
\\&=\dim_\qq H^*\big(\mathbf{B}\big( \mathbb{U}(n,\zz)/\zz^{\oplus d(n,k)} \big);\qq)\big)
\end{align*}
\end{enumerate}
\end{lemma}
\begin{prf}
Although this is not the model we use in the fibration, up to homotopy, we may suppose that
\begin{align*}
\mathbf{B} \mathbb{U}(n,\zz)=\frac{\mathbb{U}(n,\rr)}{\mathbb{U}(n,\zz)}
\end{align*}
In \cite[Theorem 1.1., Proposition 1.2]{Dwy85} it is proved that $H_*\big(\frac{\mathbb{U}(n,\rr)}{\mathbb{U}(n,\zz)};\zz\big)$ contains a free $\zz$-summand of rank $n!$, which then under $\zz_{(p)}$-localisation completely spans the localised homology due to $p\geq n-1$ (and this inequality is sharp, as there are examples with $p<n-1$ where homology contains $p$-torsion---see \cite[p.~524]{Dwy85}). By universal coefficients the same holds true for cohomology. Indeed, we obtain that
\begin{align*}
\operatorname{Ext}^1_{\zz_{(p)}}\bigg(H_*\bigg(\frac{\mathbb{U}(n,\rr)}{\mathbb{U}(n,\zz)};\zz{(p)}\bigg),\zz_{(p)}\bigg)=0
\end{align*}
as the homology is a free $\zz_{(p)}$-module.

The statements for torus cohomology in (1) are trivial and just stated for the sake of completeness in this context. Via tensoring with the flat module $\qq$ the rationalised results hold via universal coefficients.

\bigskip

It remains to derive Inequality \eqref{eqn05}, i.e.~$n!<2^{d(n,k)}$ for $n\geq 49$ and $k\geq n/2+1$. Using the Stirling formula and setting $d:=n-k\leq n/2-1$ we obtain that
\begin{align*}
n!&<\big(1+\frac{1}{11n}\big)\cdot\sqrt{2 \pi n}\cdot \big(\frac{n}{e}\big)^n<n^{n}
\intertext{for $n>1$ and}
2^{d(n,k)}&=2^{\frac{(n-k+1)(n-k+2)}{2}}=2^{\frac{d^2+3d+2}{2}}>2^{d^2/2}
\end{align*}
We solve that
\begin{align*}
2^{d^2/2}\geq n^n=2^{n\cdot \log_2 n} \iff d\geq \sqrt{2n \log_2n}
\end{align*}
For $k=\lceil n/2\rceil +1$ we obtain that $d=\lfloor n/2\rfloor-1$, and
\begin{align*}
\lfloor n/2\rfloor-1\geq  \sqrt{2n \log_2n}
\end{align*}
for $n\geq 49$. As a consequence we obtain that
\begin{align*}
\frac{\dim H^*\big(\mathbf{B} \mathbb{U}(n,\zz);\qq \big)}{\dim H^*(T^{d(n,k)} \times \rr^{c(n,k)};\qq)}\leq \frac{2^{n\cdot \log_2 n}}{2^{(\lceil n/2\rceil-1)^2/2}} \xto{n\to \infty} 0
\end{align*}
Clearly, an asymptotic formula with $\zz_{(p)}$-coefficients has to be phrased as below Theorem \ref{theoC}.
\end{prf}

\begin{proof}[\textsc{Proof of Theorem \ref{theoC}}]
The proof is a mere combination of the preceding lemmas. In Lemma \ref{lemma01} for $k\geq n/2+1$ we constructed the nilpotent fibre bundle \eqref{eqn02} of connected nilpotent rationally elliptic spaces with rational homotopy concentrated in degree $1$
\begin{align}
T^{d(n,k)} \times \rr^{c(n,k)} \hto{} \mathbf{B} \mathbb{U}(n,\zz) \to \mathbf{B}\big( \mathbb{U}(n,\zz)/\mathbb{U}_k(n,\zz) \big)
\end{align}
In Lemma \ref{lemma04}
for a prime $p\geq n-1$ we showed that the base space has finite dimensional cohomology with $R$-coefficients. Moreover, we proved in Inequalities \eqref{eqn05} and \eqref{eqn07} the estimate for the rank of the cohomologies of fibre and total space as well as the asymptotic formula.
\end{proof}

\section{Extending Theorem \ref{theoC} to highly connected examples}\label{sec04}

We shall now investigate the examples from Theorem \ref{theoC} from a purely rational point of view, which will illustrate nicely what is ``going on'' from a different, ``dual'' angle. We do so by constructing and discussing a rational minimal Sullivan model. As a next step, we then use certain alterations of the respective models in order to construct highly connected examples which have exactly the same behaviour of their rational cohomology as the examples from Theorem \ref{theoC}. That is, rationally, we can heavily extend the constructed examples.

\bigskip

 We construct a probably well-known rational model of the bundles as follows.
\begin{lemma}\label{lemma03}
A rational Sullivan model of the nilpotent fibration \eqref{eqn02} for $k\geq n/2+1$ (and, consequently, of fibre, base and total space respectively) is given by
\begin{align*}
\bigg(\Lambda \langle x_{i,j}\rangle_{i,j \leq n \atop 1\leq i-j< k-1}, \dif\bigg) \hto{} \bigg(\Lambda \langle x_{i,j}\rangle_{ i,j \leq n \atop i-j\geq 1} ,\dif\bigg) \to  \bigg(\Lambda \langle x_{i,j}\rangle_{i,j \leq n \atop i-j\geq k-1},0\bigg)
\end{align*}
with $\deg x_{i,j}=1$ for $1\leq i,j\leq n$ and the differential specified by
\begin{align*}
\dif x_{i,j}=-\sum_{j< l< i} {x_{l,j}\cdot x_{i,l}}
\end{align*}
\end{lemma}
We remark that, more generally, dropping the condition $k\geq n/2+1$, and considering the analog of the bundle \eqref{eqn02}, which clearly then has a more complicated fibre, the differential $\bar \dif$ on the model of the fibre becomes
\begin{align*}
\bar\dif x_{i,j}=-\sum_{j<l < i \atop l-j\geq k-1 , i-l\geq k-1} {x_{l,j}\cdot x_{i,l}}
\end{align*}
\begin{prf}
Due to Theorem \ref{theo05} a minimal model of the nilmanifold \linebreak[4] $\mathbb{U}(n,\rr)/\mathbb{U}(n,\zz)$ is generated by the dual of the Lie algebra of $\mathbb{U}(n,\rr)$, i.e.~it is of the form
\begin{align*}
(\Lambda \mathfrak{u}(n,\rr)^*,\dif)
\end{align*}
In the following we construct it explicitly drawing upon the description below Theorem \ref{theo05}.

Denote by $X_{i,j}=(x_{i,j})$ the matrix from $\mathfrak{u}(n,\rr)$ (with $1\leq j<i\leq n$) which has $x_{i,j}$ at position $(i,j)$ as its only non-trivial entry.

We compute that
\begin{align}\label{eqn08}
[X_{i,j},X_{s,t}]=
\begin{cases}
-X_{i,t} & \textrm{if } j=s, i\neq t\\
X_{s,j} & \textrm{if } i=t, j\neq s\\
0 &  \textrm{otherwise} \\
\end{cases}
\end{align}
The $X_{i,j}$ dualise to a basis $(x_{i,j})$ (in degree $1$ of the model). The basis $X_{i,j}$ is ordered by the order of central series extensions as depicted. We successively add the $X_{i,j}$ in the order $X_{2,1}$, $X_{3,2}$, \ldots, $X_{n,n-1}$, $X_{3,1}$, \ldots, $X_{n,n-2}$, \ldots, $X_{n-1,1}$, $X_{n,2}$, $X_{n,1}$.

Since we only consider upper triangular matrices, i.e.~$i>j$, $s>t$, and as we only deal with Lie brackets of two elements ordered in this fashion, it turns out that the first equation in \eqref{eqn08} does not apply.

Hence with this ordering and according to \eqref{eqn12} the Lie bracket transcribes to the differential, which is dual to the brackets, as
\begin{align*}
\dif x_{i,j}= -\sum_{j< l< i} x_{l,j}\cdot x_{i,l}
\end{align*}
for $i>j$.

\bigskip

The structure of $\bar \dif$ follows by projection. Hence, in the special case when $k\geq n/2+1$ this yields a trivial differential, since we derive that $\dif x_{i,j}$ can only be non-trivial if $i-j=(i-l)+(l-j)\geq 2(k-1)\geq n$. The algebra $\bigg(\Lambda \langle x_{i,j}\rangle_{i,j \leq n \atop i-j\geq k-1},0\bigg)$ then clearly is a model of the fibre torus $\mathbb{U}_k(n,\rr)/\mathbb{U}_k(n,\zz)$. From Theorem \ref{theo03} we cite that a model of the fibre is given by the fibre of the model.

Hence a model of the fibration \eqref{eqn02} is given as
\begin{align*}
(\Lambda W,\dif) \hto{} \bigg(\Lambda W\otimes \Lambda \langle x_{i,j}\rangle_{i,j \leq n \atop i-j\geq k-1},\dif\bigg)\to  \bigg(\Lambda \langle x_{i,j}\rangle_{i,j \leq n \atop i-j\geq k-1},0\bigg)
\end{align*}
for $(\Lambda W,\dif)$ a minimal model of the base  $\mathbf{B}\big( \mathbb{U}(n,\zz)/\mathbb{U}_k(n,\zz) \big)$. Since rational homotopy groups are all concentrated in degree $1$, the Sullivan model of the total space is actually minimal (not only as a relative model). We derive the isomorphism of minimal Sullivan algebras
\begin{align}\label{eqn09}
\bigg(\Lambda W\otimes \Lambda \langle x_{i,j}\rangle_{i,j \leq n \atop i-j\geq k-1},\dif\bigg)\cong \bigg(\Lambda \langle x_{i,j}\rangle_{ i,j \leq n \atop i-j\geq 1} ,\dif\bigg)
\end{align}
That is, we obtain the model of the fibre torus by taking the quotient with $(\Lambda W,\dif)$. By minimality we directly derive that
\begin{align*}
(\Lambda W,\dif)\cong \bigg(\Lambda \langle x_{i,j}\rangle_{i,j \leq n \atop 1\leq i-j< k-1}, \dif\bigg)
\end{align*}
For example, this can be obtained by first dividing out on both sides in \eqref{eqn09} the homotopy groups corresponding to the center, i.e.~$x_{n,1}$, by tensoring with $\qq[s x_{n,1}]$ ($s$ denotes the suspension, degree shift by $+1$), perturbing the differentials on $x_{n,1}$ with the summand $sx_{n,1}$, and by splitting off the contractible algebra $(\Lambda \langle x_{n,1}, sx_{n,1}\rangle, x_{n,1}\mapsto sx_{n,1})$ up to isomorphism. This induces an isomorphism on the quotient---and in a next step, we proceed in the same way with the generators of the center of the quotient, $x_{n-1,1}$ and $x_{n,2}$, and so forth. Having so divided out the fibre torus on both sides, the desired isomorphism remains. This completes the proof.
\end{prf}
Note that this generalises the well-known model of the Heisenberg group. From this model we can draw several nice observations complementing the original construction:
\begin{rem}\label{rem02}
\begin{enumerate}
\item The center of $\mathbb{U}(n,\rr)$ corresponds to the center of the Lie algebra, which itself dualises exactly to $x_{n,1}$, as this is the only element not hit by a differential (even after applying automorphisms of the model). Hence, this centre is clearly $1$-dimensional. See the proof of Proposition \ref{prop02} for the analog ``dual'' observation on the center of the Lie group $\mathbb{U}(n,\rr)$. From Theorem \ref{theo06} we derive that $\rk_0(\B \mathbb{U}(n,\zz))=1$.
\item The condition $k\geq n/2+1$ tells us exactly that the projection to the fibre yields a rational torus, i.e.~the differential $\bar \dif$ actually vanishes. For any larger subgroup the induced differential $\bar \dif$ is not trivial reflecting a non-trivial extension cocycle. That is, for example for $k=\lceil n/2\rceil$ the differential  $\bar \dif x_{n,1}$ has the non-trivial summand $-x_{n,\lfloor n/2\rfloor}\cdot x_{\lfloor n/2\rfloor,1}$.
\item Following this observation, we see that in the total space the differentials of the torus are twisted with again torus cohomology. More precisely, the constructed bundles are \emph{not principal} torus bundles. The latter come with a classifying map to $\B T$ with minimal model $(H^*(\B T),0)\xto{\simeq} \APL(\B T)$ a polynomial algebra generated in degree $2$.

    Suppose now that $\mathbf{B} \mathbb{U}(n,\zz)$ (just rationally) has the structure of a principal $T^{\frac{(n-k+1)(n-k+2)}{2}}$-bundle over $\mathbf{B}\big( \mathbb{U}(n,\zz)/\mathbb{U}_k(n,\zz))$. Then, in form of the model of the Borel fibration,
    \begin{align*}
    \mathbf{B} \mathbb{U}(n,\zz)\hto{} \mathbf{B}\big( \mathbb{U}(n,\zz)/\mathbb{U}_k(n,\zz))\to \B T^{\frac{(n-k+1)(n-k+2)}{2}}
    \end{align*}
    the space $\mathbf{B}\big( \mathbb{U}(n,\zz)/\mathbb{U}_k(n,\zz))$ has a Sullivan model of the form
    \begin{align*}
    \bigg(\Lambda \langle x_{i,j}\rangle_{ i,j \leq n \atop i-j\geq 1}\otimes \qq\big[sx_{i,j}\big]_{i,j \leq n \atop i-j\geq k-1} ,\dif\bigg)
    \end{align*}
    Due to the fact that rational homotopy groups of the total space are concentrated in degree $1$ the linear part of the differential is induced by its linear part $\dif_0 x_{i,j}=sx_{i,j}$ for $ i-j\geq k-1$. It is easy to see that this, however, infringes the property $\dif^2=0$ in all our cases (due to the depicted twisting of torus differentials). For example, the summand $x_{n,\lfloor n/2\rfloor}\cdot x_{\lfloor n/2\rfloor,1}$ of $\dif x_{n,1}$ differentiates to
    \begin{align*}
    \qquad \quad\dif \big( x_{n,\lfloor n/2\rfloor}\cdot x_{\lfloor n/2\rfloor,1}\big)=sx_{n,\lfloor n/2\rfloor}\cdot x_{\lfloor n/2\rfloor,1}-x_{n,\lfloor n/2\rfloor}\cdot sx_{\lfloor n/2\rfloor,1}\neq 0
    \end{align*}
    and this injectivity of the differential clearly still holds when taking into account all of its summands of this form. We point the reader to Proposition \ref{prop02} for a simpler proof of the non-realisability via almost-free torus actions on reasonable spaces.
\item Morally speaking, the total spaces are constructed as successive rational circle fibrations. This need for iterative steps, i.e.~the perturbed differentials, effectively thwart the growth of the rational cohomology, which turns out to be responsible for making the total space cohomology smaller than the overall fibre cohomology.
\end{enumerate}
\end{rem}

\bigskip

Next we use this model to construct further examples via spatial realisation which share the same rational cohomological properties. We may now, however, shift degrees in order to obtain simply-connected and actually arbitrarily highly connected examples. For this we adapt the model of the fibration
\begin{align*}
(\Lambda W,\dif) \hto{} \bigg(\Lambda W\otimes \Lambda \langle x_{i,j}\rangle_{i,j \leq n \atop i-j\geq k-1},\dif\bigg)\to  \bigg(\Lambda \langle x_{i,j}\rangle_{i,j \leq n \atop i-j\geq k-1},0\bigg)
\end{align*}
with $(\Lambda W,\dif)=\big(\Lambda \langle x_{i,j}\rangle_{i,j \leq n \atop 1\leq i-j< k-1}, \dif\big)$ by inductively altering degrees \emph{yet preserving differentials} as follows. Set
\begin{align*}
\deg x_{i,j}:=2 \kappa +1
\end{align*}
for $i-j=k-1=1$ on the off-diagonal with $k=2$. Suppose that $i-j+1=k$. Then set
\begin{align}\label{eqn10}
\nonumber\deg x_{i,j}:=&\deg x_{i',j'}+(2\kappa+1) -1 \textrm{ for } i'-j'+1=k-1
\\=&\deg x_{i',j'}+\deg x_{i'',j''}-1 \textrm{ whenever } (i'-j')+(i''-j'')=k-1
\end{align}
for $\kappa \geq 0$. We obtain that for $i-j+1=k$
\begin{align*}
\deg x_{i,j}= (k-1)\cdot 2\kappa+1
\end{align*}

(That is, this corresponds to starting from degree $2\kappa +1$ for the off-diagonal elements with $k=2$ and successively adding $2\kappa$ to the degrees of the respective generators corresponding to higher off-diagonals. In other words, we shift the second off-diagonal characterised by $i-j=1$ with $(2\kappa+1)$ and make this consistent with differentials by correspondingly shifting the degrees of the generators of higher \emph{lower degree}---drawing on the lower degree in rational homotopy---i.e.~on the respective successive higher off-diagonals.) Since all shifted degrees are odd, since the models we consider are coformal (i.e.~the differentials are quadratic), and since they are homogeneous with respect to lower degree, i.e.~
\begin{align*}
\dif x_{i,j}\in \langle x_{i',j'}\cdot x_{i'',j''}\rangle_{(i'-j')+(i''-j'')=i-j}
\end{align*}
by \eqref{eqn10}, this degree shift is compatible with differentials. For $\kappa=0$ we just recover the original grading---all degrees are equal to $1$ in this case.

Using spatial realisation (see \cite[Proposition 17.9, p.~248]{FHT01}) we may find a nilpotent fibration having the constructed one as its model. Hence, summarising our efforts, we have proved the highly connected rational analogue of Theorem \ref{theoC}.
\begin{theo}\label{theo02}
Suppose that $k\geq n/2+1$. Then for each $\kappa>0$ there is a nilpotent fibration of $(2\kappa)$-connected nilpotent elliptic rational spaces
\begin{align*}
F \hto{} E \to B
\end{align*}
with $H^*(F;\qq)$ the free exterior algebra on $\frac{(n-k+1)(n-k+2)}{2}$ many generators (in odd degrees), i.e.~in particular $\dim H^*(F;\qq)=2^{\frac{(n-k+1)(n-k+2)}{2}}$, and with $\dim H^*(B;\qq)<\infty$ such that
\begin{align*}
\dim H^*(E;\qq) < 2^{\frac{(n-k+1)(n-k+2)}{2}}
\end{align*}
whenever $n-k\geq \sqrt{2n \log_2n}$. In particular, both restrictions for $k$ are satisfied for $k=\lfloor n/2\rfloor +1$ and $n\geq 49$. In this case it then even holds that
\begin{align*}
\lim_{n\to \infty}\frac{\dim H^*(E;\qq)}{\dim H^*(F;\qq)}=0
\end{align*}
\end{theo}
\vproof
Clearly, on a rational space, reduced cohomology with integral coefficients agrees with reduced cohomology with rational ones.
For more on this spatial realisation see Remarks \ref{rem05} and \ref{rem01}.

\section{The rational structure of the examples in Theorem \ref{theoA}}\label{sec02}

We define a rational fibration
\begin{align}\label{eqn01}
(\Lambda \langle a,b\rangle,0)\hto{} (\Lambda \langle a,b\rangle \otimes \Lambda \langle x_1,\ldots, x_r\rangle,\dif) \to  (\Lambda \langle x_1,\ldots, x_r\rangle,\bar\dif=0)
\end{align}
with $\deg a=\deg b=\deg x_i=1$ for all $1\leq i\leq r$. The differential is well-defined by
\begin{align*}
\dif a&=\dif b=0,\\
\dif x_1&=a\cdot b,\\
\dif x_2&=a\cdot x_1, \ldots,\\
\dif x_i&=a \cdot x_{i-1} \textrm{ for $2\leq i\leq k$}
\end{align*}
By construction, this constitutes a $\Lambda$-extension in the notation of \cite[p.~xix]{FHT15} or what is usually called a \emph{rational fibration}. We denote its  total space by $X_r$.

\begin{rem}\label{rem05}
Before we apply a different method of passing from algebra to geometry in Section \ref{sec05}, which is better suited and stronger in our setting, let us here quickly address the usual way of realising rational structures geometrically. In particular, this serves as the blueprint for constructing highly connected examples---see Remark \ref{rem01}.

Hence we apply spatial realisation to this rational fibration in order to obtain a genuine one. That is, from \cite[p.~115]{FHT15} we cite that
\begin{align*}
|\langle (\Lambda \langle a,b\rangle,0)|\rangle \leftarrow |\langle (\Lambda \langle a,b\rangle \otimes \Lambda \langle x_1,\ldots, x_r\rangle,\dif) |\rangle\leftarrow  |\langle(\Lambda \langle x_1,\ldots, x_r\rangle,\bar\dif=0)\rangle|
\end{align*}
is a fibre bundle of CW-complexes and a Serre fibration fitting into the following commutative diagram (see \cite[Proposition 3.14, p.~115]{FHT15})
\begin{align*}
\tiny
\xymatrix{
\APL(|\langle (\Lambda \langle a,b\rangle,0)|) \ar[r]
& \APL( |\langle (\Lambda \langle a,b\rangle \otimes \Lambda \langle x_1,\ldots, x_r\rangle,\dif) |\rangle) \ar[r]
& \APL(|\langle(\Lambda \langle x_1,\ldots, x_r\rangle, 0)\rangle|)\\
(\Lambda \langle a,b\rangle,0)  \ar[u]
\ar@{^{(}->}[r]& (\Lambda \langle a,b\rangle \otimes \Lambda \langle x_1,\ldots, x_r\rangle,\dif)  \ar[u]
\ar[r]
& (\Lambda \langle x_1,\ldots, x_r\rangle, 0) \ar[u]
}
\end{align*}

Since the vector spaces underlying our models are all concentrated in degree $1$ and are finite dimensional, it follows from \cite[Theorem 5.4, p.~166]{FHT15} that the vertical morphisms are quasi-isomorphisms, i.e.~Sullivan models which are actually even minimal by degree.

Let us further see that the constructed fibration is quasi-nilpotent and nilpotent, actually. For this we observe that all of fibre, base and total space are nilpotent spaces. This follows from \cite[Remark, p.~15]{FH17}. Indeed, clearly, for each space say with model $(\Lambda W^1,\dif)$ there is some $N$ such that $W^1_N=W^1$ where $W_0^1\In \ldots \In W_r^1\In \ldots$ is the filtration defined by $W_0^1=W^1\cap \ker \delta$, $W_{r+1}=\delta^{-1}(W^1\cap W_r^1)$, and $\delta w$ is the component of $\dif w$ lying in $W^1\wedge W$---which in our case just yields $\delta=\dif$. From Remark \ref{rem03}.(4) we cite that due to the nilpotency of base and total space the whole fibration is nilpotent, which is, in particular, known to be quasi-nilpotent (see Remark \ref{rem03}.(2)).

The minimal model of the base space, $(\Lambda \langle a,b\rangle,0)$, is the one of a $2$-torus, and we may use the continuous map $T^2 \to |\langle (\Lambda \langle a,b\rangle,0)\rangle |$, which can be considered as the map from $T^2$ to its rationalisation $T^2_\qq$ (see \cite[p.~108, Theorem 17.12(ii), p.~254]{FHT01}), in order to pull back the fibre bundle. Due to \cite[Proposition 4.4.3, p.~69]{MP12} the total space $E$ is nilpotent, and this yields a nilpotent bundle
\begin{align*}
 |\langle(\Lambda \langle x_1,\ldots, x_r\rangle,\bar\dif=0)\rangle|\hto{} E \to T^2
\end{align*}
of nilpotent spaces.
\end{rem}

\begin{prop}\label{prop01}
The total dimension of the cohomology $\dim H^*(X_r)$ is as depicted in Table \ref{table01}.

\setlength{\tabcolsep}{10pt}
\renewcommand{\arraystretch}{1}
\begin{table}[h]
\centering \caption{total dimension $\dim H^*(X_r;\qq)$}
\label{table01}
\begin{tabular}{c  |  c | c}
$r$ & $2^r$ & $\dim H^*(X_r;\qq)$\\
\hline
$0$ & $1$&  $3$\\
$1$ & $2$ & $6$\\
$2$ & $4$ & $8$\\
$3$ & $8$ & $12$\\
$4$ & $16$ & $16$\\
\hdashline
$5$ & $32$ & $26$\\
$6$ & $64$ & $40$\\
$7$ & $128$ & $64$\\
$8$ & $256$ & $104$\\
$9$ & $512$ & $180$
\end{tabular}
\end{table}
\end{prop}
\begin{prf}
This follows from computations using \textsc{SageMath 9.0}. For the convenience of the reader we provide the vector space generators of $H^*(X_5;\qq)$, corresponding to the minimal $r=5$ for which $\dim H(X_r)<2^r$. These are
\begin{align*}
\{&[1],[a], [b], [bx_1], [x_1x_2 - bx_3], [bx_1x_2],[x_2x_3 - x_1x_4 + bx_5], [ax_5], [ax_2x_3],
\\&
 [bx_2x_3 - bx_1x_4], [x_1x_2x_3 - bx_2x_4 + bx_1x_5], [ax_3x_4],[bx_1x_2x_3],[ax_4x_5],
\\&
 [bx_1x_3x_4 - bx_1x_2x_5], [ax_2x_3x_4], [x_1x_2x_3x_4 - bx_2x_3x_5 + bx_1x_4x_5],
\\&
[ax_2x_3x_5], [ax_1x_2x_3x_4], [bx_1x_2x_3x_4], [ax_3x_4x_5], [ax_1x_2x_4x_5], [ax_2x_3x_4x_5],
\\&
 [ax_1x_2x_3x_4x_5], [bx_1x_2x_3x_4x_5], [abx_1x_2x_3x_4x_5]\}
\end{align*}
\end{prf}

\begin{rem}\label{rem01}
By attributing different degrees to the generators of $X_r$ in \eqref{eqn01}, for example like $\deg a=\deg b=3$, $\deg x_i=5+2\cdot(i-1)$, one may produce simply-connected (or even arbitrarily highly connected) examples similar to the ones in Theorem \ref{theoA}. Spatial realisation of this fibration can be done along the lines of Remark \ref{rem05}. Note, however, that for simply-connected spaces spatial realisation becomes much more simple (see \cite[Chapter 17, p.~237]{FHT01}).
Hence one can construct nilpotent fibre bundles of the form
\begin{align*}
 |\langle(\Lambda \langle x_1,\ldots, x_r\rangle,0)\rangle|\hto{} E \to B
\end{align*}
where the cohomology of $B$ is freely generated by two odd-degree generators, and the $x_i$ are concentrated in odd degrees. One may choose the cohomology of fibre base and total space to be arbitrarily highly connected, and clearly $\dim H^*(E;\qq)<2^r$ for $5\leq r\leq 9$.
\end{rem}

\section{Finishing the proofs of Theorems \ref{theoCC} and \ref{theoA}}\label{sec05}

In Sections \ref{sec04} and \ref{sec02}, in particular, we understood the rational structures of the examples considered in Theorems \ref{theoC} and \ref{theoA}.  Based upon the rational structure of the examples from Theorem \ref{theoC} we use the geometric realisation depicted below in order to construct the examples from Theorem \ref{theoCC}. As for the constructions in Theorem \ref{theoA} we build upon the rational examples from Section \ref{sec02} and use the subsequent geometric realisation as well. This allows us to finish the proofs of the main theorems in this Section.

In fact, we have seen that rationally the models we considered result as the total spaces of a sequence of rational $\s^1$-fibrations. The base space of these iterative fibrations is, as we have seen in Section \ref{sec01} (or actually by applying the subsequent arguments just in order to construct a suitable base space) a smooth manifold with rational homotopy concentrated in degree $1$. For example, we may decompose the fibration \eqref{eqn01} underlying the examples from Theorem \ref{theoA} into a sequence of rational $\s^1$-fibrations
\begin{align*}
(\Lambda \langle a,b\rangle \otimes \Lambda \langle x_1,\ldots, x_{i-1}\rangle,\dif) \hto{} (\Lambda \langle a,b\rangle \otimes \Lambda \langle x_1,\ldots, x_i\rangle,\dif) \to  (\Lambda \langle x_i\rangle,0)
\end{align*}
and our iterative construction starts over a genuine $2$-torus with model freely generated by $a,b$.

This is a very convenient setting in order to apply a different classical method of geometric realisation. In all cases
we shall construct a finite length sequence (in the last case $r$ many) of genuine principal $\s^1$-fibrations $\s^1\hto{}E_i\to E_{i-1}$ over the designated base $E_0$ (in the example $E_0=T^2$) such that in each step the model of the total space, $E_i$, of the geometric construction is the total space of the corresponding rational circle fibration. (Again, in the example, the total space of the $i$-th-construction step has total space rationally equivalent to $(\Lambda \langle a,b\rangle \otimes \Lambda \langle x_1,\ldots, x_i\rangle,\dif)$.)

For this, we start by recalling that principal $\s^1$-bundles over a CW-complex $X$ are classified by homotopy classes of maps $[X,\B \s^1]$ with $\B \s^1=\cc\pp^\infty=K(\zz,2)$ (and $\E \s^1=\s^\infty$). Since cohomology is described by the Eilenberg--MacLane spectrum, this allows for an easier description. That is, since $H^2(X;\zz)\cong \operatorname{Pic}(X)$ is isomorphic to the Picard group of complex line bundles (together with the tensor product, and the isomorphism given by attributing first Chern classes), we hence may choose a classifying map $\phi\co X\to \B \s^1$ with $H^2(\cc\pp^\infty;\zz)=\langle z\rangle$ and the property, without restriction, that $H^*(\phi)(z)$ is an integral representative of the differential in the Sullivan model, i.e.~of the Euler class in the rational fibration. (Again, in the first step of the example, we have that $H^*(T^2;\zz)=H^*(T^2;\zz)=\Lambda \langle a,b\rangle$---once having chosen $a$, $b$ to be generating integral representatives of the corresponding rational classes above---and we set $H^*(\phi)(z)=a\cdot b$, since $\dif x_1=a\cdot b$.) First, by cellular approximation this map may be supposed to map to a finite dimensional complex projective space, by smooth approximation it then may be assumed to be smooth. Now use $\phi^*$ to pull back a complex line bundle $L$ to $X$.

We finally pass to the associated smooth $\s^1$-bundle $\mathbb{S}L$, which then, by construction, has Euler class $\chi(\mathbb{S}L)$ a multiple of the rational Euler class. (In the showcase $\chi(\mathbb{S}L)=a\cdot b$.)
From \cite[Example 4, p.~202]{FHT01} it follows that the $i$-th rational circle fibration is a model of the $i$-th geometrically constructed one; in particular, the total space is a model (even minimal by degree) of the geometric total space. (Our adaptations of choosing multiples do not change the isomorphism type of the underlying Sullivan models.)

Hence, iterating this construction, starting with $\s^1\hto{} E_1\to E_0$ and consequently progressing via bundles $\s^1\hto{} E_i\to E_{i-1}$ leaves us with a geometric realisation of the rational fibration in the form $F\hto{} E\to B$, where $B=E_0$ and the fibre $F$ results from the successive circle fibres. Indeed, we obtain that
\begin{align*}
E_i&=E_{i-1} \times_{\cc\pp^\infty} \s^\infty
\\&=(E_{i-2} \times_{\cc\pp^\infty} \s^\infty) \times_{\cc\pp^\infty} \s^\infty
\\&=(\ldots (E_{0} \times_{\cc\pp^\infty} \s^\infty) \times_{\cc\pp^\infty} \s^\infty) \times_{\cc\pp^\infty} \ldots \times_{\cc\pp^\infty} \s^\infty )\ldots )\times_{\cc\pp^\infty} \s^\infty
\end{align*}
i.e.~successive bundles $F_i\hto{} E_i\to E_0$ with fibres
\begin{align*}
F_i=(\ldots (\{\textrm{pt}\} \times_{\cc\pp^\infty} \s^\infty) \times_{\cc\pp^\infty} \s^\infty) \times_{\cc\pp^\infty} \ldots \times_{\cc\pp^\infty}\s^\infty )\ldots )\times_{\cc\pp^\infty} \s^\infty
\end{align*}
After the last step (say the $j$-th one) in this finite iteration this stabilises and yields the bundles $F:=F_j\hto{} E:=E_j\to B=E_0$. The rational model tells us that all these $F_i$, and $F$ in particular, have the rational type of a genuine torus. It remains to see that $F$ is a diffeomorphism torus.

Since we are working in the smooth category, our spaces are smooth manifolds and our maps are smooth as well. The fibre $F$ is a smooth manifold and is constructed as an iterative principal circle bundle. From \cite[Theorem 1]{Bel18} we cite
\begin{theo}\label{theo08}
A manifold is an iterated principal circle bundle if and only if it is a
nilmanifold.
\end{theo}
From \cite[Remark 3.21.3, p.~121]{FOT08} we recall that a nilmanifold of the rational type of a torus is diffeomorphic to a torus; hence so is $F=T$. (Alternatively, not drawing on these two results, we may see that also integrally the respective first Chern classes of the pullback bundles vanish resulting in a torus fibre.)

Observe further, that we may start constructing the base space $B$ as an iterated principal circle bundle as well. (By the theorem this is already the case once we start with a nilmanifold $B$.) Hence so becomes the total space $E$. By the theorem, both spaces then are nilmanifolds. In particular, they are nilpotent spaces and the fibration is nilpotent as well (see Remark \ref{rem03}.4). This nilpotence of $E$ (and hence of the fibration in view of the remark and the nilpotence of $B$) also just follows from \cite[Theorem 1.1, p.~8]{Hil76},
since $\pi_*(E)=\pi_1(E)$, and since the fundamental group is constructed via a sequence of central extensions over a nilpotent group.

\begin{proof}[\textsc{Proof of Theorem \ref{theoCC}}]
The proof goes completely along the lines of the one of Theorem \ref{theoC}. We just use the smooth nilpotent fibre bundle we constructed in this section instead of the one considered there. The estimates \eqref{eqn05} (with rational coefficients) and \eqref{eqn07} clearly still apply, since the considered fibrations are rationally equivalent.
\end{proof}

\begin{proof}[\textsc{Proof of Theorem \ref{theoA}}]
Again we use the geometric constructions from this section applied to the rational constructions from Section \ref{sec02}, i.e.~to the rational fibration \eqref{eqn01}. The rational estimates are provided by Proposition \ref{prop01}. This proves the first part of the theorem. As for the second one, we already observed in the introduction that taking product bundles is enough to cover all fibre torus ranks from $5$ onwards.
\end{proof}


\section{Non-realisability via almost free torus actions}\label{sec03}

In this section we show that none of the presented examples are realisable by almost free torus actions. However, we first need to make sense out of this statement.

\begin{conv}\label{conv01}
For any of the considered fibrations from Theorems \ref{theoC}/\ref{theoCC}, \ref{theo02}, \ref{theoA}
and Remark \ref{rem01} we show that they are ``not realisable by an almost free torus action''. By this we mean two-fold.
\begin{enumerate}
\item In the cases of Theorems \ref{theoC}/\ref{theoCC} and \ref{theoA} the concrete models of the rational torus fibrations cannot be the rational models of torus fibrations corresponding to almost free torus actions. (This is not restricted to the total spaces being reasonable.)
\item More generally: any reasonable space (in the definition and distinction from Section \ref{subsec03}) with the rational homotopy type of the total space in a fibration given in one of Theorems \ref{theoC}/\ref{theoCC}, \ref{theo02}, \ref{theoA}, and Remark \ref{rem01} does not admit an almost free action of the corresponding torus.
\end{enumerate}
\end{conv}

\begin{rem}\label{rem04}
It is important to note that although our examples have homotopy tori as fibres, it follows from our discussion below the algebraic toral rank conjecture \ref{conj03} that already the ``correct rational data'' usually can be realised by a smooth action on a closed manifold. Hence the following results showing non-realisability are essential in order to understand that the torus bundles are not ``almost free torus actions in disguise'', and hence cannot be counter-examples to the toral rank conjecture.
\end{rem}

\subsection{The ``matrix examples''}\label{subsec01}

We deal with the examples of Theorem \ref{theoC}/\ref{theoCC} and Theorem \ref{theo02}. In Remark \ref{rem02}.(3) we already showed the first aspect of Convention \ref{conv01}, namely that a concrete fibration from Theorem \ref{theoC}/\ref{theoCC} cannot come from an almost free torus action.

In Remark \ref{rem02}.(1) we moreover already read off from the structure of the minimal model that $\rk_0(\B \mathbb{U}(n,\zz))=1$, i.e.~that the rational toral rank of any example from Theorem \ref{theoC} is one. Hence the largest almost free torus action on a reasonable space of the rational type of the total space can come from an $\s^1$, consequently. This, however, required the understanding of the rational model in Lemma \ref{lemma03}. For the convenience of the reader let us provide a simpler proof here.
\begin{prop}\label{prop02}
The examples $\mathbf{B} \mathbb{U}(n,\zz)$ presented in Theorem \ref{theoC} satisfy
\begin{align*}
1=\rk_0  \mathbf{B} \mathbb{U}(n,\zz)=C(\pi_*(\Omega  \mathbf{B} \mathbb{U}(n,\zz))\otimes \qq,[ \cdot, \cdot])
\end{align*}
The same holds for the rational torus ranks and the centers of the homotopy Lie algebras of the total spaces from Theorem \ref{theo02}.
\end{prop}
\begin{prf}
Up to homotopy (which clearly has no effect on the considered rational invariants) we may assume that $\mathbf{B} \mathbb{U}(n,\zz)= \mathbb{U}(n,\rr)/\mathbb{U}(n,\zz)$. The latter is a nilmanifold with associated Lie algebra $\mathfrak{u}(n,\rr)$; due to Theorem \ref{theo06} we know that the rational toral rank equals the center, i.e.~
\begin{align*}
\rk_0(\mathbb{U}(n,\rr)/\mathbb{U}(n,\zz))= \dim C(\mathfrak{u}(n,\rr))
\end{align*}
The Lie algebra of the center of the group is the center of the Lie algebra. It is easy to check that $C(\mathbb{U}(n,\rr))$ consists exactly of those matrices with all off-diagonal entries $0$ except for possibly the entry at position $(n,1)$. That is,
\begin{align*}
\rk_0(\mathbb{U}(n,\rr)/\mathbb{U}(n,\zz))= \dim C(\mathfrak{u}(n,\rr))= \dim C(\mathbb{U}(n,\rr))=1
\end{align*}
By the remark following Theorem \ref{theo07} this agrees with the dimension of the center of the homotopy Lie algebra.
Note that all the relevant examples in Theorem \ref{theoC} were constructed for tori of rank much larger than $1$.

\bigskip

The examples from Theorem \ref{theo02} result from positive degree shifts out of the ones from Theorem \ref{theoC}. The rational homotopy groups of their total spaces $E$ remain concentrated in odd degrees, and their minimal models are (non-gradedly) isomorphic to the one of $\mathbf{B} \mathbb{U}(n,\zz)$. As a consequence we derive that the center of their homotopy Lie algebra remains $1$-dimensional. Having realised the total spaces as finite CW-complexes (see \cite[p.~91]{FOT08}) and then using Theorem \ref{theo07} we conclude that also their rational toral rank is at most one, and then actually equal to one, as above. (Indeed, we can just easily specify a rational fibration with base $\B\s^1$ satisfying $H^2(\B \s^1;\qq)=\langle t\rangle$ and with the element corresponding to the suspension of this central element, say $sx_{n,1}\in \pi_*(E)\otimes \qq$, perturbed by adding to it a suitable power of $t$ under the differential. The discussion below Conjecture \ref{conj03} yields the geometric realisation.)
\end{prf}
As a consequence, it follows that the differences between the ranks of the torus bundles we construct and the actual rational toral rank (which is just $1$) becomes arbitrarily large. An orientable $\s^1$-bundle is realisable by an $\s^1$-action, and the toral rank conjecture trivially holds in this case.

\subsection{The ``abstract examples''}\label{subsec02}

In this section we deal with the examples from Theorem \ref{theoA} respectively from Remark \ref{rem01}.

Let us first compute the toral rank of these examples $E$ given as the specific total spaces $T^k\hto{} E\to T^2$ as in the first part of Theorem \ref{theoA}. (The examples from Remark \ref{rem01} are (non-gradedly) isomorphic whence the analog arguments apply leading to the identical computations and results.)

We do so by actually computing the centres of their homotopy Lie algebras $L:=(\pi_*(\Omega E)\otimes \qq,[\cdot, \cdot ])$. We then use that we realised them as finite CW-complexes, respectively that, abstractly, there is a finite CW-complex realising the spaces $E$ (see \cite[p.~91]{FOT08}) such that we are in the realm of reasonable spaces. Then we invoke Theorem \ref{theo07} in order to derive the estimate $\rk_0 E\leq \dim C(\pi_*(\Omega E) \otimes \qq,[\cdot, \cdot ])$, as the rational homotopy of $E$ is clearly concentrated in odd degree equal to $1$.
\begin{prop}\label{prop03}
It holds that $1=\rk_0 E= \dim C(\pi_*(\Omega E) \otimes \qq)$.
\end{prop}
\begin{prf}
Recall their minimal model
\begin{align}
(\Lambda \langle a,b\rangle,0)\hto{} (\Lambda \langle a,b\rangle \otimes \Lambda \langle x_1,\ldots, x_r\rangle,\dif) \to  (\Lambda \langle x_1,\ldots, x_r\rangle,\bar\dif=0)
\end{align}
with $\deg a=\deg b=\deg x_i=1$ for all $1\leq i\leq r$, and
\begin{align*}
\dif a&=\dif b=0,\\
\dif x_1&=a\cdot b,\\
\dif x_2&=a\cdot x_1, \ldots,\\
\dif x_i&=a \cdot x_{i-1} \textrm{ for $2\leq i\leq r$}
\end{align*}
We want to determine that the only element in the center is the rational homotopy group corresponding to $s^{-1}x_r$, the desuspension (degree shift by $-1$) of the top degree generator. That is, we show that this is the only element (up to multiples) commuting with all others. Our model is clearly coformal, i.e.~all differentials are quadratic. Since the quadratic part of the differential in the minimal Sullivan model is dual to the Lie bracket (see \cite[p.~46]{FHT15}), we just have to show that the only element not hit by a differential is $x_r$. Since
\begin{align*}
\dif x_i=a x_{i-1} \textrm{ for $2\leq i\leq r$}\qquad \textrm{and} \qquad \dif x_1=ab
\end{align*}
we derive that for $x\in \langle a,b, x_1,\ldots, x_{r-1}\rangle$ we have that $ax \in \im \dif$, and, up to duality, $[s^{-1} a, s^{-1} x_{i-1}]=-s^{-1}x_{i}\neq 0$ for $i\leq r$. Moreover, clearly, $s^{-1}x_r$ lies in the center. As in the proof of Proposition \ref{prop02}, $x_r$ can be used to construct a free $\s^1$-fibration. The result follows.
\end{prf}

\begin{rem}
As we have seen in Propositions \ref{prop02} and \ref{prop03} all our examples have rational toral rank equal to the dimension of the centre of the homotopy Lie algebra equal to $1$. Hence confirming the toral rank conjecture in these cases is a complete triviality. Nonetheless, let us make the following observation: All the examples are coformal and satisfy that the centre of the homotopy Lie algebra is concentrated in its highest degree. In this situation (irrespective of the actual dimension of the centre) the toral rank conjecture follows from \cite[Proposition 3.2, Remark 3.3, p.~8, p.~9]{DS88}.
\end{rem}

\bigskip

This finishes Aspect (2) of Convention \ref{conv01}. Before we provide several further non-realisability results and interpretations let us compute Aspect (1), the concrete non-realisability for the given fibration.

So let us make a non-realisability argument fine-tuned to the rational structure of the fibration---similar to the one given in Remark \ref{rem02}.(3) for the ``matrix examples''. So let us see that $X_r$, the total space of the rational fibration \eqref{eqn01}, cannot even rationally be the fibre of a classifying fibration of a principal $T^r$-bundle with base rationally a $T^2$, i.e.~of a hypothetical rational fibration
\begin{align*}
X_r\hto{} (X_r)_{T^r}\simeq_\qq T^2\to \B T^r
\end{align*}
i.e.~that it does not fit into the setting of the algebraic toral rank conjecture (see Conjecture \ref{conj03}) for $r>1$. Since $\B T^r$ is formal with cohomology $H^*(\B T^r)=\qq[t_1, \ldots, t_r]$, it follows that any twisted differential of a model of $(X_r)_T$ by degree has to be of the form
\begin{align*}
\dif a&=\dif b=0 \\
\dif x_1&=ab+\tilde t_1\\
\vdots&\\
\dif x_i&=ax_{i-1}+ \tilde t_i
\end{align*}
for $2\leq i\leq r$ and with $\tilde t_i\in \langle t_1,\ldots, t_r\rangle$ of degree $2$. (The fact that $\dif a=\dif b=0$ follows from a short computation very similar to the following one.) For $i\geq 3$ we compute
\begin{align*}
\dif^2 x_i=\dif(ax_{i-1}+\tilde t_i)=-a(ax_{i-2}+\tilde t_{i-1})=-a \tilde t_{i-1}
\end{align*}
since $a^2=0$. We deduce that $\tilde t_{i-1}=0$; the analogue computation shows that $\tilde t_1=0$. For dimension reasons, it follows that $(X_r)_{T^r}$ has non-trivial rational homotopy groups in degree $2$ contradicting the fact that it is rationally a $T^2$.

(More specifically, for reasonable spaces, using Theorem \ref{theo04} one may argue that, without restriction, $\qq[s \tilde t_{i-1}]$ injects into $(X_r)_T$ contradicting finite dimensional cohomology unless $\tilde t_{i-1}=0$.)

\bigskip

Finally, we remark that, clearly, the rational toral rank of a product is at least as large as the sum over the rational toral ranks of the factors (with strict inequality actually possible---see \cite[Example 3.3, p.~204]{JL04}). Hence the examples constructed for the second part of Theorem \ref{theoA} certainly have rational toral rank larger than just $1$.

\subsection{Further interpretation and context}\label{subsec04}
Let us finally put the results into further context.
In \cite[Theorem A, p.~3]{AZ19} we proved (using the $H^*(\B T)$-module structure) that the toral rank conjecture holds for compact Hausdorff spaces $X$ whenever the Borel construction $X_T$ is formal. So again, we realised the examples from Theorem \ref{theoA} by finite CW-complexes and assume they admit an almost free action by $T$. Then the orbit space is necessarily a $T^2$, and the formality of $T^2$ gives a contradiction, i.e.~it yields the non-realisability of the torus bundles via almost free torus actions.

\bigskip

Yet another alternative and maybe more enlightening way of seeing non-realisability is the following. Again, but only exemplarily, we discuss the examples from Theorem \ref{theoA}, i.e.~the rational fibration \eqref{eqn01}. In the case of realisability, $T\hto{} X_r\to (X_r)_T$ is the pullback fibration of the universal fibration $T\hto{} \E T\to \B T$ (again letting the cohomology of the base be generated by the $t_i$ in degree $2$). The space $\E T$ is contractible, and we have that $\dif_u s^{-1}t_i=t_i$ for $1\leq i\leq r$ in the model of this universal fibration with differential $\dif_u$. The differential on $T\hto{} X_r\to (X_r)_T\simeq_\qq T^2$ is induced from the universal one by composition with the induced map of the classifying morphism $\phi\co |\langle (X_r)_T\rangle|\to \B T$, i.e.~$\dif=H^*(\phi)\circ \dif_u$. Since $(X_r)_T\simeq_\qq T^2$ and hence $\dim H^2((X_r)_T)=1$, the map $H^2(\phi)$ has a large kernel (given a large torus $T$). Hence so has $\dif$, already on $\langle x_1,\ldots, x_r\rangle$, which clearly is a contradiction.

This shows the following: The given examples owe their existence to the fact that the classifying space $\B \aut_0(T_\qq)\not\simeq_\qq \B T_\qq$ of the considered rational torus fibration over $X_r$ (thus respectively for an arbitrary $T$-fibration) is necessarily rationally more complicated than the classifying space of a principal $T$-bundle. (Recall that $\aut_0(T_\qq)$ denotes the group-like space of self-maps of $T_\qq$ homotopic to the identity, and $\B \aut_0(T_\qq)$ is its classifying space.) Or, in other words, the classifying maps for the fibrations we consider, $T^2_\qq\to \B\aut_0(T_\qq)$ do not factor over maps $T^2_\qq\to (\B T)_\qq\to \B\aut_0(T_\qq)$.

\bigskip

Recall further that (for example see \cite[Theorem 2.2, p.~354]{LS20}) for a \emph{simply-connected} space $X$ of finite type the algebra $\Der(\Lambda V,\dif)$ is a Lie model for $\B \aut_0(X_\qq)$ (where $(\Lambda V,\dif)$ is a minimal model of $X$, $X_\qq$ its rationalisation, and $\Der(\Lambda V,\dif)$ is the differential graded Lie algebra of negative degree derivations of $(\Lambda V,\dif)$). In particular, this implies that there is an isomorphism of graded Lie algebras
\begin{align*}
\pi_*(\Omega \B \aut_0(X_\qq))\cong H_*(\Der(\Lambda V,\dif))
\end{align*}
(see \cite[Theorem 2.1, p.~354]{LS20})---after restricting to connected Lie algebras only.

\bigskip

We remark that using this description and replacing an $\s^1$-torus action by an $\Sp(1)\cong \SU(2)\cong \s^3$-torus action---denoted by $T_{\s^3}$---we do obtain that
\begin{align*}
\B \aut_0((T_{\s^3})_{\qq})\simeq_\qq \B T_{\s^3}
\end{align*}
Indeed, in this case we see that the derivations $\delta$ of degree $-3$ on a minimal model $(\Lambda y_1,\ldots, y_k\rangle,0)$ of $T_{\s^3}$ given by $\delta(y_i)=1$ and $\delta y_j=0$ for $i\neq j$ generate exactly $(H^*(\B T_{\s^3}),0)=(\qq[sy_1,\ldots, sy_k],0)$, a model for $\B T_{\s^3}$.

Note however, that for different simply-connected Lie group fibres over simply-connected base spaces there are again fibrations which are not principal. In \cite[Example 4, p.~220]{FHT01} there can be found such an example with fibre $\SU(2)\times \SU(3)$. Analyzing it, this phenomenon is clearly due to the (non-trivial) derivation of degree $-2$ induced by sending the top degree algebra generator of $H^*(\SU(3);\qq)$ to the one of $H^*(\SU(2);\qq)$, clearly yielding a rational homotopy group not in $\pi_*(\Omega\B\SU(2))\otimes \qq\oplus \pi_*(\Omega\B\SU(3))\otimes \qq$, and showing that respective classifying spaces are rationally different.

Also from this perspective another relation to the Halperin conjecture---stating that a fibration of simply-connected spaces with positively elliptic fibre $F$ is totally non-homologous to zero---becomes visible. The Halperin conjecture holds if and only if $\pi_*(\B \aut_0(F))=\pi_\even(\B \aut_0(F))$ (see \cite[Theorem A]{Mei82}).


\def\cprime{$'$}

\pagebreak \

\vfill

\begin{center}
\noindent
\begin{minipage}{\linewidth}
\small \noindent \textsc
{Manuel Amann} \\
\textsc{Institut f\"ur Mathematik}\\
\textsc{Differentialgeometrie}\\
\textsc{Universit\"at Augsburg}\\
\textsc{Universit\"atsstra\ss{}e 14 }\\
\textsc{86159 Augsburg}\\
\textsc{Germany}\\
[1ex]
\footnotesize
\textsf{manuel.amann@math.uni-augsburg.de}\\
\textsf{www.uni-augsburg.de/de/fakultaet/mntf/math/prof/diff/team/dr-habil-manuel-amann/}
\end{minipage}
\end{center}

\end{document}